\documentclass[11pt]{amsart}

\parskip=\smallskipamount 

\usepackage{amsfonts} 
\usepackage{amsmath} 
\usepackage{amsthm} 
\usepackage{amssymb} 
\usepackage{pst-plot} 

\newcommand{\set}[2]{\{#1;\ #2\}}
\renewcommand{\d}{\partial}
\newcommand{\dbar}{\bar{\partial}}

\DeclareMathOperator{\supp}{supp}

\newcommand{\NN}{\mathbb{N}}
\newcommand{\ZZ}{\mathbb{Z}}
\newcommand{\RR}{\mathbb{R}}
\newcommand{\CC}{\mathbb{C}}
\newcommand{\OO}{\mathcal{O}}
\renewcommand{\AA}{\mathcal{A}}

\newtheorem{theorem}{Theorem}[section]
\newtheorem{proposition}[theorem]{Proposition}
\newtheorem{corollary}[theorem]{Corollary}
\newtheorem{lemma}[theorem]{Lemma}

\theoremstyle{definition}
\newtheorem{definition}[theorem]{Definition}

\newtheorem{remark}[theorem]{Remark}

\numberwithin{equation}{section}

\begin{document}

\title{Approximation of holomorphic mappings on 1-convex domains}

\author{Kris Stopar}
\address{Institute of Mathematics, Physics and Mechanics, University of Ljubljana, Jadranska ulica 19, 1000 Ljubljana, Slovenia}
\email{kris.stopar@fmf.uni-lj.si}


\begin{abstract}
Let $\pi :Z\to X$ be a holomorphic submersion of a complex manifold $Z$ onto a complex manifold $X$ and $D\Subset X$ a 1-convex domain with strongly pseudoconvex boundary. We prove that under certain conditions there always exists a spray of $\pi $-sections over $\bar{D}$ which has prescribed core, it fixes the exceptional set $E$ of $D$, and is dominating on $\bar{D}\setminus E$. Each section in this spray is of class $\mathcal{C}^k(\bar{D})$ and holomorphic on $D$. As a consequence we obtain several approximation results for $\pi $-sections. In particular, we prove that $\pi $-sections which are of class $\mathcal{C}^k(\bar{D})$ and holomorphic on $D$ can be approximated in the $\mathcal{C}^k(\bar{D})$ topology by $\pi $-sections that are holomorphic in open neighborhoods of $\bar{D}$. Under additional assumptions on the submersion we also get approximation by global holomorphic $\pi $-sections and the Oka principle over 1-convex manifolds. We include an application to the construction of proper holomorphic maps of 1-convex domains into $q$-convex manifolds.
\end{abstract}

\maketitle

\section{Introduction}
\label{sec:introduction}

Modern methods in complex analysis are in large part based on the solutions of the inhomogeneous Cauchy-Riemann equations with estimates. Among the most important ones are the $L^2$-method developed by H\"orman\-der, Kohn and others, and the method of integral operators with holomorphic kernels developed by Henkin, Arrelano, Lieb, \O vrelid and others. These methods are linear and have to be modified for the use in nonlinear problems. One of such nonlinear problems is the classical Oka-Grauert theory on the classification of holomorphic $G$-bundles over Stein spaces \cite{Gra}, \cite{Gra1}. The linearization of the basic problem in this theory was achieved by approximation theorems for holomorphic mappings into Lie groups, by Cartan splitting lemma, and by the gluing method for such maps. This theory was generalized significantly in the work of Gromov \cite{Gro} to the case where Lie groups are replaced with a much bigger class of manifolds, the so called elliptic manifolds, i.e., manifolds with a global dominating spray. It was in his work that the theory of approximation and gluing of holomorphic objects using sprays became a more general theory for the first time.

In the recent years the theory of sprays and their applications has been developed even further through the works of Forstneri\v c, Drinovec-Drnov\v sek, Prezelj and others. A part of this theory is based on local sprays of sections (which may exist in a more general setting) rather than global sprays (which exist only on certain manifolds). Local sprays have proven to be a very useful tool when dealing with all kinds of problems, such as:
\begin{itemize}
\item Oka-Grauert-Gromov theory on Stein manifolds with strongly pseudoconvex boundary \cite{Dri&For1},
\item approximation by holomorphic mappings in a neighborhood of the domain \cite{Dri&For1},
\item construction of open Stein neighborhoods of embedded Stein subvarieties with strongly pseudoconvex boundary \cite{For1},
\item describing the Banach manifold structure of spaces of holomorphic mappings \cite{For1},
\item new constructions of proper holomorphic mappings \cite{Dri&For2},
\item new approaches to the treatment of disc functionals \cite{Dri&For3},
\end{itemize}
and many others.

In this paper we extend the technique of sprays, developed in \cite{Dri&For}, \cite{Dri&For1}, \cite{Dri&For2} and \cite{For}, to obtain some of the above results in the more general context of 1-convex manifolds. Our results extend the work of Henkin and Leiterer \cite{Hen&Lei} and Prezelj \cite{Pre}. One of the highlights is the Oka principle for sections of holomorphic fiber bundles with Oka fibers over compact 1-convex domains with strongly pseudoconvex boundary (cf.\ Theorem \ref{opr}).

The structure of the paper is the following. In Sec.\ \ref{sec:sprays} we recall the notion of a (dominating) holomorphic spray of sections, and we introduce a certain Condition $\mathcal{E}$ concerning the existence of sprays of sections with a given core on a neighborhood of the exceptional set of a 1-convex manifold. The main results of the paper are described in Sec.\ \ref{sec:main}. In Sec.\ \ref{sec:preliminaries} we prove some technical lemmas which allow us to construct a bounded linear solution operator for the inhomogeneous $\dbar$-equation for certain $(0,1)$-forms on 1-convex domains with strongly pseudoconvex boundary (Theorem \ref{bso}). In Sec.\ \ref{sec:Cartan} we define 1-convex Cartan pairs and prove the generalized Cartan splitting lemma for such pairs with $\mathcal{C}^k$-estimates up to the boundary (Theorem \ref{st}). In Sec.\ \ref{sec:gluing} we show how this splitting lemma can be used to glue sprays of sections over 1-convex Cartan pairs (Proposition \ref{gs}). This gluing method is the main tool for the rest of the article. In Sec.\ \ref{sec:existence} we construct sprays of sections over 1-convex domains. In Sec.\ \ref{sec:approximation} we prove a local Mergelyan-type approximation result for sections of holomorphic submersions onto 1-convex domains (Theorem \ref{la}, see also Theorem \ref{lam}) and a version of the Oka principle with approximation for sections of certain holomorphic fiber bundles on 1-convex domains which extend smoothly to the boundary (Theorem \ref{opr}).

\section{Sprays of sections over 1-convex domains}
\label{sec:sprays}

Let $X$ be a complex manifold. An open relatively compact subset $D$ of $X$ is said to be a \emph{1-convex domain with strongly pseudoconvex boundary} of class $\mathcal{C}^k$, $k\ge 2$, if there exists an open neighborhood $W\subset X$ of the boundary $\d D$ and a strictly plurisubharmonic function $\varphi :W\to \RR $ of class $\mathcal{C}^k$ such that $D\cap W=\set{x\in W}{\varphi (x)<0}$. If $W$ can be chosen to be a neighborhood of the whole set $\bar{D}$, then $D$ is a Stein domain with strongly pseudoconvex boundary. Similarly, a manifold $X$ is called \emph{1-convex} if there exists a smooth exhaustion function $\varphi :X\to \RR $ which is strictly plurisubharmonic outside a compact subset of $X$.

By a classical characterization of 1-convexity due to Narasimhan \cite{Nar1}, a 1-convex domain (or manifold) $D$ is a proper modification of a Stein space at a finite set of points; i.e., there exist a Stein space $Y$, a proper holomorphic surjection $R:D\to Y$ with connected fibers, and a finite set $F\subset Y$ such that $Y\setminus F$ is nonsingular, $R:D\setminus \pi^{-1}(F)\to Y\setminus F$ is a biholomorphism, and $R_*(\OO _D)=\OO _Y$. The map $R:D\to Y$ (or just the space $Y$) is called the \emph{Remmert reduction} of $D$, and its degeneracy set $E:=\set{x\in D}{\dim_x R^{-1}(R(x))>0}$ is called the \emph{exceptional set} of $D$.

Suppose that $Z$ is a complex manifold. Let $k\in \{0,1,2,\ldots,\infty\}$. We shall say that a map $f:\bar{D}\to Z$ is of class $\AA ^k$ if it is of class $\mathcal{C}^k(\bar{D})$ and holomorphic on $D$. We will consider sections of a surjective map $\pi :Z\to \bar{D}$ which is either the restriction to $\bar{D}$ of a holomorphic submersion $\pi :\tilde{Z}\to X$ or a fiber bundle of class $\AA ^k$. Let us recall the definition of such a bundle.

\begin{definition}
Let $X$ and $Y$ be complex manifolds, and let $D\subset X$ be a domain with $\mathcal{C}^l$ boundary for some $l\ge 2$. A fiber bundle $\pi :Z\to \bar{D}$ with fiber $Y$ is said to be of class $\AA^k$ $(k\in\mathbb{Z}_+)$ if every point $z\in \bar{D}$ admits a relatively open neighborhood $U$ in $\bar{D}$ and a fiber bundle isomorphism $\phi :Z|_{U}\to U\times Y$ which is of class $\mathcal{C}^k$ and holomorphic over $U\cap D$.
\end{definition}

Let $Z_x:=\pi ^{-1}(x)$ denote the fiber over $x\in \bar{D}$. For each $z\in Z$ denote by $VT_zZ:=\ker D\pi (z)$ the tangent space to the fiber $Z_{\pi (z)}$ at $z$, also called the \emph{vertical tangent space} of $Z$ at $z$. The {\em vertical tangent bundle} $VTZ$ with fibers $VT_z Z$ is a holomorphic vector subbundle of the tangent bundle $TZ$.

We recall the notion of a holomorphic spray of sections; cf.\ Def.\ 5.9.1 in \cite[p.\ 215]{For3}.

\begin{definition}\label{s}
Let $D\subset X$ be a domain with $\mathcal{C}^l$ boundary ($l\ge 2$) and $0\in P\subset \CC ^N$ an open set. A (local) \emph{spray of $\pi $-sections of class $\AA^k(D)$} is a $\mathcal{C}^k$ map $f:\bar{D}\times P\to Z$ which is holomorphic on $D\times P$ and satisfies
$$\pi (f(x,t))=x, \quad x\in \bar D,\ t\in P.$$
Such $f$ is \emph{dominating} on $K\subset \bar{D}$ if the partial differential
$$\d _t|_{t=0}f(x,t):T_0\CC ^N\cong \CC ^N\longrightarrow VT_{f(x,0)}Z$$
is surjective for every $x\in K$. We say that $f$ \emph{fixes} a subset $L\subset \bar{D}$ if
$$f(x,t)=f(x,0), \quad x\in L,\ t\in P.$$
The $\pi$-section $f_0:=f(\cdot ,0)$ is called the \emph{core} of the spray of $\pi $-sections $f$. If $P=\CC ^N$ we say that $f$ is a {\em global spray of $\pi $-sections of class $\AA ^k(D)$}.
\end{definition}

Every global spray of sections is of course also a local spray of sections; hence we shall omit the word `local' when talking about sprays of sections that may or may not be global.

Similarly one defines a holomorphic spray of sections over an open domain in the base space. The domination condition will be very important in the process of gluing sprays of sections (see Proposition \ref{gs}). In the case of a canonical projection $\pi :X\times Y\to X$ we can identify $\pi $-sections with maps from $X$ to $Y$ and sprays of $\pi $-sections with sprays of maps.

Sprays of sections should not be confused with the classical notion of a spray due to Gromov \cite{Gro} which we now recall. Let $p:V\to Z$ be a holomorphic vector bundle. Denote by $V_z:=p^{-1}(z)$ the fiber over a point $z\in Z$ and by $0_z$ the zero element in $V_z$. Notice that $V_z$ is in a natural way a $\CC $-linear subspace of the tangent space $T_{0_z}V$.

\begin{definition}\label{sg}
Let $\pi :Z\to X$ be a holomorphic submersion and $D$ a domain in $X$. A (local) \emph{$\pi $-spray} on $Z|_D=\pi ^{-1}(D)$ is a quadruple $(V,\Omega ,p,s)$, where $p:V\to Z|_D$ is a holomorphic vector bundle, $\Omega \subset V$ is an open set containing the zero section, and $s:\Omega \to Z|_D$ is a holomorphic map satisfying
$$s(V_z)\subset Z_{\pi (z)}, \quad  s(0_z)=z, \quad z\in Z|_D.$$
Such a spray is \emph{dominating} on $K\subset Z|_D$ if the restriction of the differential
$$Ds(0_z):V_z\to VT_zZ$$ 
is surjective for every $z\in K$. If $\Omega =V$ we say that the triple $(V,p,s)$ is a global $\pi $-spray on $Z|_D$.
\end{definition}

There exists a dominating spray of sections with a prescribed core over any relatively compact Stein subset of $X$ \cite[Lemma 5.10.4, p.\ 220]{For3} (see also Proposition \ref{eosp} below and the preceding discussion). If $X$ is 1-convex, one could expect the existence of sprays which are dominating over the complement of the exceptional set $E\subset X$. An attempt in this direction was made in \cite[\S 4]{Pre} (see also \cite[Corollary 2.6]{Pre&Sla}). However, the proof in \cite{Pre} is incomplete (see Sec.\ \ref{sec:existence} below or the Appendix in \cite{For&Lar}). The main problem is to find such sprays in a small neighborhood of $E$. To avoid this difficulty we introduce the following condition.

\medskip \noindent
\textbf{Condition $\boldsymbol{\mathcal{E}}$.} Let $X$ be a 1-convex manifold with the exceptional set $E$. We say that a holomorphic submersion $\pi :Z\to X$ satisfies \emph{Condition $\mathcal{E}$} if for every open neighborhood $U_0\subset X$ of $E$ and every holomorphic $\pi $-section $f_0:U_0\to Z$ there exists an open neighborhood $U\subset U_0$ of $E$ and a holomorphic spray of $\pi $-sections $f:U\times P\to Z$ with core $f_0|_U$ which is dominating on $U\setminus E$ and fixes $E$.
\medskip

We expect that Condition $\mathcal{E}$ is always satisfied, but we are currently unable to prove it. The most important sufficient conditions implying Condition $\mathcal{E}$ are given by the following proposition (for the proof see  Sec.\ \ref{sec:existence}).

\begin{proposition}\label{ce}
Let $X$ be a 1-convex manifold with the exceptional set $E$.
\begin{enumerate}
\item Let $\pi :Z\to X$ be a holomorphic submersion. Suppose that there exists an open neighborhood $D\subset X$ of $E$ and a $\pi $-spray $(V,\Omega ,p,s)$ on $Z|_D$ which is dominating on $Z|_{D\setminus E}$. Then $\pi $ satisfies Condition~$\mathcal{E}$.
\item If $Y$ is an elliptic manifold in the sense of Gromov \cite{Gro} then the canonical projection $\pi :X\times Y\to X$ satisfies Condition $\mathcal{E}$.
\item If $Y$ is a complex manifold such that global holomorphic vector fields on $Y$ span the tangent space $T_y Y$ at each point $y\in Y$, then the canonical projection $\pi :X\times Y\to X$ satisfies Condition $\mathcal{E}$.
\end{enumerate}
\end{proposition}

Assuming Condition $\mathcal{E}$ we can find sprays of sections over larger domains as in the following theorem.

\begin{theorem}\label{eost}
Let $D$ be a relatively compact 1-convex domain with strongly pseudoconvex boundary of class $\mathcal{C}^l$ ($l\ge 2$) and exceptional set $E$ in a complex manifold $X$, and let $\pi:Z\to \bar{D}$ be either a fiber bundle of class $\AA ^k$ ($k\le l$) or the restriction to $\bar D$ of a holomorphic submersion $\tilde{Z}\to X$. Suppose that the submersion $\pi :Z|_D\to D$ satisfies Condition $\mathcal{E}$. Given a $\pi $-section $f_0:\bar{D}\to Z$ of class $\AA ^k$ there exists a spray of $\pi $-sections $f:\bar{D}\times P\to Z$ of class $\AA ^k$ with core $f_0$ which is dominating on $\bar{D}\setminus E$ and fixes $E$.
\end{theorem}

Theorem \ref{eost} is proved in Sec.\ \ref{sec:existence}.

The only thing that is incorrect or incomplete in \cite{Pre} is the proof of the existence of a holomorphic spray of $\pi $-sections, which has prescribed core and is dominating on the complement of the exceptional set (see Sec.\ \ref{sec:existence} below for further explanation). The existence of such sprays is now given by Theorem \ref{eost}. Thus, under the additional assumption that the submersion $\pi $ satisfies Condition $\mathcal{E}$, the results of [24] remain valid.

\section{The main results}
\label{sec:main}
We now describe the main results of the paper on approximation of maps over 1-convex domains. We begin with a local Mergelyan-type approximation theorem. In the case when $D$ is a strongly pseudoconvex Stein domain (i.e., with a trivial exceptional set), this is Theorem 1.2 in \cite{Dri&For1}.

\begin{theorem}\label{lam}
Let $D$ be a relatively compact 1-convex domain with strongly pseudoconvex boundary of class $\mathcal{C}^l$ ($l\ge 2$) and exceptional set $E$ in a complex manifold $X$ and let $Y$ be a complex manifold such that the canonical projection $X\times Y\to X$ satisfies Condition $\mathcal{E}$. Every map $f:\bar{D}\to Y$ of class $\AA ^k$ ($k\le l$) can be approximated in $\mathcal{C}^k(\bar{D})$ by maps which are holomorphic in open neighborhoods of $\bar{D}$ and agree with $f$ on $E$.
\end{theorem}

In Sec.\ \ref{sec:approximation} we extend Theorem \ref{lam} to sections of certain fiber bundles over 1-convex domains (see Theorem \ref{la}).

For an approximation result by global maps we need to make additional assumptions on the target manifold which allows sufficient flexibility for holomorphic maps into it. The most significant such condition is the Oka property. The simplest of several equivalent definitions of an {\em Oka manifold} is the following (cf.\ Definition 5.4.1 in \cite[p.\ 192]{For3}).

\begin{definition}\label{om}
A complex manifold $Y$ is an \emph{Oka manifold} if any holomorphic map from an open neighborhood of a compact convex set $K\subset \CC ^m$ ($m\in \NN $) to $Y$ can be approximated uniformly on $K$ by entire maps $\CC ^m\to Y$.
\end{definition}

This property is also referred to as the {\em convex approximation property}, or CAP for short \cite[p.\ 192]{For3}. It implies the full Oka principle for maps of Stein manifolds to $Y$ (see Theorem 5.4.4 in \cite[p.\ 193]{For3}). 

Our next result is the Oka principle with approximation for sections of fiber bundles of class $\AA^k$ with Oka fibers over compact 1-convex domains. The analogous result in the Stein case was obtained in \cite[Theorem 1.7]{Dri&For1}. The classical case of holomorphic fiber bundles with homogenous fibers over 1-convex manifolds is due to Henkin and Leiterer \cite{Hen&Lei}. (For fiber bundles over Stein manifolds see Grauert \cite{Gra1} and Grauert and Kerner \cite{Gra&Ker}.)

\begin{theorem}\label{opr}
Let $X$ be a complex manifold and let $D\Subset D'\Subset X$ be bounded 1-convex domains with strongly pseudoconvex boundary of class $\mathcal{C}^l$ ($l\ge 2$) and the same exceptional set $E$. Let $\pi :Z\to \bar{D}'$ be a fiber bundle of class $\AA ^k$ $(k\in \{0,1,\ldots,l\})$ with Oka fiber such that the submersion $\pi :Z|_D\to D$ satisfies Condition $\mathcal{E}$. If $\bar{D}$ is $\OO(D')$-convex, then for every continuous $\pi $-section $f_0:\bar{D}'\to Z$ which is of class $\AA^k$ on $\bar{D}$ there exists a homotopy $f_t:\bar{D}'\to Z$ of continuous $\pi $-sections of class $\AA ^k$ on $\bar{D}$ such that $f_t$ is close to $f_0$ in $\mathcal{C}^k(\bar{D})$ and agrees with $f_0$ on $E$ for every $t\in [0,1]$, and the section $f_1$ is of class $\AA ^k$ on $\bar{D}'$.
\end{theorem}

Theorem \ref{opr} still holds if we replace $\bar{D}$ by an arbitrary compact $\OO (D')$-convex set $K\subset D'$ and assume that the initial section $f_0$ is holomorphic in a neighborhood of $K$. By using Theorem \ref{opr} inductively, we immediately obtain the corresponding Oka principle in the open case, that is, for sections of a holomorphic fiber bundle $Z\to X$ with Oka fiber over a 1-convex base $X$ (see Corollary \ref{opc} and compare with the main result in \cite{Pre}).

A complex manifold $Y$ of complex dimension $n$ is said to be $q$-convex if it admits a $\mathcal{C}^2$ exhaustion function $\varphi :Y\to \RR $ that is $q$-convex outside some compact set $K\subset Y$ (i.e., the Levi form of $\varphi$ has at least $n-q+1$ positive eigenvalues at each point of $Y\setminus  K$). If $K$ can be chosen empty then the manifold $Y$ is said to be $q$-complete. Our next theorem is another application of the technique of gluing holomorphic sprays of maps.

\begin{theorem}\label{phm}
Let $X$ and $Y$ be complex manifolds of complex dimensions $m$ and $n$, respectively, and let $D\Subset X$ be a 1-convex domain with smooth strongly pseudoconvex boundary. Assume that the canonical projection $D\times Y\to D$ satisfies Condition $\mathcal{E}$, that $2m\le n$, and let $q\in \{1,\ldots ,n-2m+1\}$. Then the following hold:
\begin{enumerate}
\item If $Y$ is $q$-convex then there exists a proper holomorphic map $D\to Y$.
\item If $Y$ is $q$-complete then every continuous map $\bar{D}\to Y$ that is holomorphic in $D$ can be approximated, uniformly on compacts in $D$, by proper holomorphic maps $D\to Y$.
\end{enumerate}
\end{theorem}

The Stein case of Theorem \ref{phm} can be found in \cite{Dri&For2} along with a more general version and additional corollaries. By using the techniques developed in this paper, the proof in \cite{Dri&For2} easily caries over to our 1-convex case.

\section{Preliminaries}
\label{sec:preliminaries}

The following two propositions are stated in \cite[Lemma 2.2 and Satz 3.2]{Ric} in a slightly weaker form.

\begin{proposition}\label{ppf}
Let $A$ be an analytic subset of $D \subset \CC ^N$ and $K\subset A$ a compact subset possessing an open Stein neighborhood in $D$. Then there exist an open neighborhood $\Omega \subset D$ of $K$ and a smooth plurisubharmonic function $\psi :\Omega \to \RR $ with the following properties:
\begin{enumerate}
\item $\psi $ is strictly plurisubharmonic on $\Omega \setminus A$,
\item the Levi form of $\psi $ is positive definite on $A\cap \Omega $ in the directions which are not Zariski tangent to $A$, and 
\item $\psi =0$ on $A\cap \Omega$ and $\psi >0$ on $\Omega \setminus A$.
\end{enumerate}
\end{proposition}

\proof
We follow the proof of \cite[Lemma 2.2]{Ric}. Let $\mathcal{I}\subset \OO (D)$ be the sheaf of ideals defined by $A$. Choose an open Stein neighborhood $U\subset D$ of the compact set $K$ (such exists by assumption). For every $N$-tuple $f=(f_1,\ldots ,f_N)\in H^0(U,\mathcal{I})^N$ we denote by $J(f;z)$ the complex Jacobian of $f$ at the point $z\in U$. For every point $w\in U\setminus A$ there exists $f\in H^0(U,\mathcal{I})^N$ with $\det J(f;w)\neq 0$. Indeed, by Cartan's Theorem A there exists $f_0\in H^0(U,\mathcal{I})$ with $f_0(w)\neq 0$. We set $f_j(z):=\frac{(z_j-w_j)f_0(z)}{f_0(w)}$ for $z\in U$, $j=1,\ldots ,N$. Then $f=(f_1,\ldots ,f_N)$ has the desired property, since $J(f;w)=Id$. The analytic set $B=\set{z\in U}{\det J(f;z)=0 \textrm{ for all } f\in H^0(U,\mathcal{I})^N}$ is contained in $A\cap U$. Since $K$ is compact, there exist an open neighborhood $\Omega \subset D$ of $K$ and finitely many $f^{(1)},\ldots ,f^{(q)}\in H^0(U,\mathcal{I})^N$ such that $B\cap \Omega =\set{z\in U}{\det J(f^{(j)};z)=0 \textrm{ for all } j=1,\ldots ,q}$. Set $f=(f_1,\ldots ,f_m):=(f^{(1)},\ldots ,f^{(q)})$. Then the matrix $J(f;z)$ has maximal rank for all points $z\in \Omega \setminus A$. Since this does not change if we add a finite number of elements from $H^0(U,\mathcal{I})$ into $f=(f_1,\ldots ,f_m)$ we may assume that $A\cap \Omega =\set{z\in U}{f_j(z)=0 \textrm{ for all } j=1,\ldots , m}$. For $z\in \Omega $ we set
$$\psi (z):=\sum _{j=1}^mf_j(z)\overline{f_j(z)}.$$
Obviously $\psi $ is smooth, $\psi =0$ on $A\cap \Omega$, and $\psi >0$ on $\Omega \setminus A$. It is easy to see that the Levi form of $\psi $ at a point $z\in \Omega $ evaluated on a vector $v\in\CC ^N$ equals $|J(f;z)v|^2$, hence $\psi $ is plurisubharmonic on $\Omega $. Moreover, $\psi $ is strictly plurisubharmonic on $\Omega \setminus A$ since $J(f;z)$ has maximal rank there, and its Levi form is positive definite on $A\cap \Omega $ in the directions which are not Zariski tangent to $A$.
\endproof

\begin{proposition}\label{espf}
Let $A$ be an analytic subset of a domain $D \subset \CC ^N$, $K\subset A$ be a compact subset, and $\phi :A\to \RR $ be a strictly plurisubharmonic function of class $\mathcal{C}^k$ ($k\ge 2$). Then $\phi $ can be extended to a $\mathcal{C}^k$ strictly plurisubharmonic function on an open neighborhood $\Omega \subset D$ of the set $K$. \end{proposition}

\proof
The proof goes exactly as in \cite[Satz 3.2.b]{Ric}. It is clear from the proof that the degree of smoothness of the extension depends on the degree of smoothness of the starting function. 
\endproof

The following lemma gives a smooth analogue for the maximum function, which has been used by many authors (see for example \cite[\S 4]{Hen&Lei1}). We adapt it slightly for our purposes. We shall need it for the construction of a certain Stein domain (see Proposition \ref{esd}).

\begin{lemma}\label{sm}
Let $\epsilon >0$ and choose a convex function $A_{\epsilon }:\RR \to \RR $ of class $\mathcal{C}^{\infty }$ such that $A_{\epsilon }(x)=x$ for $x\ge 0$ and $A_{\epsilon }(x)=-x-\epsilon $ for $x\le -\epsilon $. Then the function $M_{\epsilon }:\RR ^2\to \RR $, defined by $M_{\epsilon }(x,y)=\frac{x+y}{2}+A_{\epsilon }(\frac{x-y}{2})$, has the following properties:
\begin{enumerate}
\item $M_{\epsilon }(x,y)=x$ for $x\ge y$,
\item $M_{\epsilon }(x,y)=y-\epsilon $ for $x\le y-2\epsilon $,
\item $M_{\epsilon }(x,y)\ge x$ and $M_{\epsilon }(x,y)\ge y-\epsilon $,
\item $0\le \max \{x,y\}-M_{\epsilon }(x,y)\le \epsilon $, and 
\item if $X$ is a complex manifold and $\Phi ,\Psi :X\to \RR $ are strictly plurisubharmonic functions of class $\mathcal{C}^k$ ($k\ge 2$) then $M_{\epsilon }(\Phi ,\Psi )$ is also a strictly plurisubharmonic function of class $\mathcal{C}^k$. 
\end{enumerate}
\end{lemma}

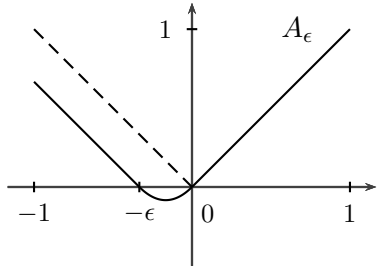
\begin{figure}[h]

\psset{unit=0.7cm}
\begin{pspicture}(-3.6,-1.6)(3.6,3.6)

\psaxes[linecolor=darkgray,labels=none,ticks=none]{->}(0,0)(-3.5,-1.5)(3.5,3.5)

\psplot[plotstyle=curve,plotpoints=100]{-3}{-1}{-1 x sub}
\psplot[plotstyle=curve,plotpoints=100]{-1}{0}{x 1 add x mul}
\psplot[plotstyle=curve,plotpoints=100]{0}{3}{x}
\psplot[plotstyle=curve,plotpoints=100,linestyle=dashed]{-3}{0}{0 x sub}

\psline(-1,-0.1)(-1,0.1)
\rput(-1,-0.5){\small{$-\epsilon $}}
\psline(-3,-0.1)(-3,0.1)
\rput(-3,-0.5){\small{$-1$}}
\psline(3,-0.1)(3,0.1)
\rput(3,-0.5){\small{$1$}}
\psline(-0.1,3)(0.1,3)
\rput(-0.5,3){\small{$1$}}
\rput(0.3,-0.5){\small{$0$}}
\rput(2,3){$A_{\epsilon }$}

\end{pspicture}

\caption{The function $A_{\epsilon }$.}
\label{figA}
\end{figure}

\begin{figure}[h]

\psset{unit=0.7cm}
\begin{pspicture}(-1.1,-2.1)(9.2,4.1)

\pspolygon[fillstyle=solid,fillcolor=lightgray,linestyle=none](4,-1.64)(0.96,-0.8)(0,0)(0.48,2.19)(1.56,3)(4,3.64)

\parametricplot[plotstyle=curve,fillstyle=solid,fillcolor=lightgray,linestyle=none]{-1.236}{-0.4}{t t 2 sub mul t 0.4 sub}
\psecurve[fillstyle=solid,fillcolor=lightgray,linestyle=none](1.92,-0.8)(0.96,-0.8)(0,0)(0.96,0.8)
\parametricplot[plotstyle=curve,fillstyle=solid,fillcolor=lightgray,linestyle=none]{0}{2.2}{0.1 t 2 exp mul t}
\psecurve[fillstyle=solid,fillcolor=lightgray,linestyle=none](0.6,1.4)(0.48,2.19)(1.56,3)(2.6,3)
\parametricplot[plotstyle=curve,fillstyle=solid,fillcolor=lightgray,linestyle=none]{2.6}{3.236}{t t 2 sub mul t 0.4 add}

\pspolygon[fillstyle=hlines,linestyle=none,hatchcolor=gray](4,4)(1.6,4)(0.4,-2)(4,-2)
\parametricplot[plotstyle=curve,fillstyle=hlines,linestyle=dashed,hatchcolor=gray]{-2}{4}{0.1 t 2 exp mul t}
\parametricplot[plotstyle=curve,fillstyle=vlines,linestyle=dashed,hatchcolor=gray]{-1.236}{3.236}{t t 2 sub mul t}

\parametricplot[plotstyle=curve]{-1.236}{-0.4}{t t 2 sub mul t 0.4 sub}
\psecurve(1.92,-0.8)(0.96,-0.8)(0,0)(0.96,0.8)
\parametricplot[plotstyle=curve]{0}{2.2}{0.1 t 2 exp mul t}
\psecurve(0.6,1.4)(0.48,2.19)(1.56,3)(2.6,3)
\parametricplot[plotstyle=curve]{2.6}{3.236}{t t 2 sub mul t 0.4 add}

\psframe[fillstyle=hlines,linestyle=none,hatchcolor=gray](5,3)(6,4)
\rput[l](6.5,3.5){$\{\Phi <0\}$}

\psframe[fillstyle=vlines,linestyle=none,hatchcolor=gray](5,1)(6,2)
\rput[l](6.5,1.5){$\{\Psi <0\}$}

\psframe[fillstyle=solid,fillcolor=lightgray,linestyle=none](5,-1)(6,0)
\rput[l](6.5,-0.5){$\{M_{\epsilon }(\Phi ,\Psi )<0\}$}

\end{pspicture}

\caption{The zero level set of the function $M_{\epsilon }(\Phi ,\Psi )$.}
\label{figM}
\end{figure}
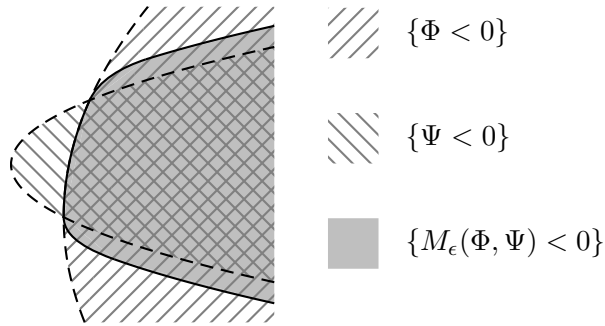

Notice that $M_{\epsilon }(x,y)$ is not symmetric in $x$ and $y$. We may think of it as a smooth maximum of $x$ and $y-\epsilon $. Figure \ref{figA} shows the graph of the function $A_{\epsilon }$ and Figure \ref{figM} shows the zero level set of the function $M_{\epsilon }(\Phi ,\Psi )$ with respect to the zero level sets of functions $\Phi $ and $\Psi $. 

\proof
The properties (i) and (ii) clearly follow from the definitions of $A_{\epsilon }$ and $M_{\epsilon }$. Notice that $\max \{x,y\}=\frac{x+y}{2}+|\frac{x-y}{2}|$. The convexity of $A_{\epsilon }$ ensures that $A_{\epsilon }(x)\ge x$ and $A_{\epsilon }(x)\ge -x-\epsilon $, which implies property (iii). Moreover, we have $0\le |x|-A_{\epsilon }(x)\le \epsilon $, which implies property (iv). For the proof of property (v) set $N=(\Phi ,\Psi )$. It is clear that $M_{\epsilon }(\Phi ,\Psi )$ is of class $\mathcal{C}^k$. Its Levi form equals
$$L_zM_{\epsilon }(\Phi ,\Psi )(v)=\langle H_{N(z)}M_{\epsilon }\cdot \d _zN(v),\d _zN(v)\rangle +D_{N(z)}M_{\epsilon }\cdot L_zN(v),$$
where $D$ and $H$ denote the differential and the Hessian operators, respectively, and $\langle \cdot ,\cdot \rangle $ denotes the hermitian inner product. We have
$$H_{N(z)}M_{\epsilon }=\frac{1}{4}A_{\epsilon }''\Big(\frac{\Phi (z)-\Psi (z)}{2}\Big)\left[\begin{array}{cc} 1 & -1 \\ -1 & 1 \end{array}\right].$$
Since $A_{\epsilon }$ is convex, the hessian $H_{N(z)}M_{\epsilon }$ is nonnegative definite. Therefore $\langle H_{N(z)}M_{\epsilon }\cdot \d _zN(v),\d _zN(v)\rangle $ is nonnegative. Furthermore, we have
$$D_{N(z)}M_{\epsilon }=\frac{1}{2}\left(1+A_{\epsilon }'\Big(\frac{\Phi (z)-\Psi (z)}{2}\Big),1-A_{\epsilon }'\Big(\frac{\Phi (z)-\Psi (z)}{2}\Big)\right).$$
The convexity of $A_{\epsilon }$ ensures that $-1\le A_{\epsilon }'\le 1$. Therefore $D_{N(z)}M_{\epsilon }$ has nonnegative entries that are not all $0$. Since $L_zN(v)=(L_z\Phi (v),L_z\Psi (v))^T$ has positive entries for $v\ne 0$, it follows that $D_{N(z)}M_{\epsilon }\cdot L_zN(v)$ is positive for $v\ne 0$, which implies property (v).
\endproof

We now show how to `extend' a Stein domain $G$ in a Stein subvariety $Y\subset \CC ^N$ to a domain in the ambient space $\CC ^N$ with certain properties. This will be helpful when dealing with 1-convex domains, since their Remmert reduction is such a domain $G$.

\begin{proposition}\label{esd}
Let $Y\subset \CC ^N$ be a Stein subvariety, $G\Subset Y$ a Stein domain with strongly pseudoconvex boundary $\d G\subset Y_{reg}$ of class $\mathcal{C}^k$ ($k\ge 2$) and $V\subset G\cap Y_{reg}$ an open set. Suppose that there exists an open neighborhood $\Omega _0'\subset \CC ^N$ of $\bar{V}$ and a holomorphic retraction $r:\Omega _0'\to Y\cap \Omega _0'$. Let $\Omega _0\Subset \Omega _0'$ be another open neighborhood of $\bar{V}$ and let $\Omega _1\subset \CC ^N$ be an open neighborhood of $\bar{G}\setminus \Omega _0$. Then there exists a Stein domain $\Omega \Subset \CC ^N$ with strongly pseudoconvex boundary of class $\mathcal{C}^k$ such that $r(\Omega \cap \Omega _0)\subset G\subset \Omega \Subset \Omega _0\cup \Omega _1$.
\end{proposition}

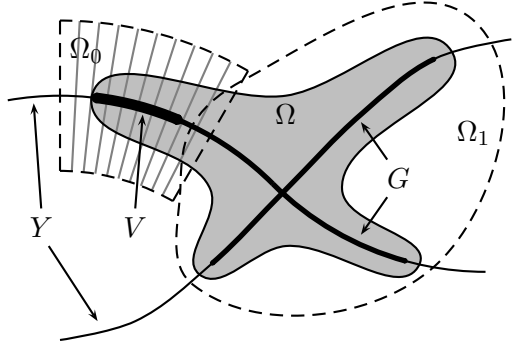
\begin{figure}[h]

\psset{unit=0.5cm}
\begin{pspicture}(-1.5,1.4)(12.5,11)
\SpecialCoor

\psecurve[fillstyle=solid,fillcolor=lightgray,linestyle=none](3,9)(8;84)(4,5.5)(3.7,3.25)(6,4)(9.55,3.5)(7.5,5.5)(10.4,9.2)(6,8)(8;84)(3,7)

\multido{\i=63+3}{9}{\psline[linecolor=gray](6;\i)(10;\i)}
\psline[linestyle=dashed](6;60)(10;60)
\psline[linestyle=dashed](6;90)(10;90)
\psarc[linestyle=dashed](0,0){6}{60}{90}
\psarc[linestyle=dashed](0,0){10}{60}{90}

\psarc(0,0){8}{45}{100}
\rput(16;45){\psarc(0,0){8}{225}{270}}
\pscurve(0,1.5)(2,2)(4,3.5)(6,5.5)(8,7.5)(10,9)(12,9.5)

\psarc[linewidth=2pt,arrows=-cc](0,0){8}{45}{84}
\rput(16;45){\psarc[linewidth=2pt,arrows=-cc](0,0){8}{225}{255}}
\psecurve[linewidth=2pt,arrows=cc-cc](2,2)(4,3.5)(6,5.5)(8,7.5)(10,9)(12,9.5)

\psarc[linewidth=4pt,arrows=cc-cc](0,0){8}{66}{84}

\psecurve(3,9)(8;84)(4,5.5)(3.7,3.25)(6,4)(9.55,3.5)(7.5,5.5)(10.4,9.2)(6,8)(8;84)(3,7)

\psccurve[linestyle=dashed](8;64)(6.4;60)(3.5,3)(10,3.5)(11,10)(9.7;60)

\rput(-0.5,4.5){$Y$}
\psline[arrows=->](-0.5,5)(-0.7,7.8)
\psline[arrows=->](-0.3,4)(1,2)
\rput(9,5.8){$G$}
\psline[arrows=->](8.7,6.3)(8,7.3)
\psline[arrows=->](8.7,5.3)(8,4.3)
\rput(2,4.5){$V$}
\psline[arrows=->](2,5)(2.2,7.5)
\rput(0.7,9.2){$\Omega _0$}
\rput(11,7){$\Omega _1$}
\rput(6,7.5){$\Omega $}

\end{pspicture}

\caption{The domain $\Omega $.}
\label{figO}
\end{figure}

\proof
Since $G\subset Y$ is a Stein domain with strongly pseudoconvex boundary of class $\mathcal{C}^k$, there exists a strictly plurisubharmonic function $\phi $ of class $\mathcal{C}^k$ defined in an open neighborhood $W'\subset Y$ of the set $\bar{G}$ (with isolated critical points) such that $G=\set{z\in W'}{\phi (z)<0}$. We may assume that $W'=(\Omega _0'\cup \Omega _1)\cap Y$.

By Proposition \ref{espf} the function $\phi $ can be extended to a strictly plurisubharmonic function of class $\mathcal{C}^k$ on some open neighborhood $\Omega _2\subset \CC ^N$ of the set $\bar{G}$. By Proposition \ref{ppf} there exist an open neighborhood $\Omega _3\subset \CC ^N$ of the set $\bar{G}$ and a plurisubharmonic function $\psi :\Omega _3\to \RR $ of class $\mathcal{C}^{\infty }$ which is strictly plurisubharmonic everywhere except on $\Omega _3\cap Y$ in the directions Zariski tangent to $Y$, such that $\psi =0$ on $Y\cap \Omega _3$ and $\psi >0$ on $\Omega _3\setminus Y$. We may assume that $\Omega _2=\Omega _3=\Omega_0'\cup \Omega _1$.

We define functions $\tilde{\Phi }=\phi \circ r+\psi $ on $\Omega _0'$ and $\Psi =\phi +M\psi $ on $\Omega _0'\cup \Omega _1$, where the constant $M$ will be chosen positive. The function $\Psi $ is clearly strictly plurisubharmonic on $\Omega _0'\cup \Omega _1$. The same holds true for the function $\tilde{\Phi }$ on $\Omega _0'$. Indeed, the functions $\phi \circ r$ and $\psi $ are plurisubharmonic there, whereby $\psi $ is strictly plurisubharmonic everywhere except on $Y$ in the Zariski tangent directions, where $\phi \circ r=\phi $ is strictly plurisubharmonic.

We would now like to glue the functions $\tilde{\Phi }$ and $\Psi $ together using Lemma \ref{sm}. To ensure a smooth transition we need to add a small correction to the function $\tilde{\Phi }$ as follows. Choose open sets $\Omega _0\Subset\Omega _0'''\Subset \Omega _0''\Subset \Omega _0'$ and a smooth function $\chi :\CC ^N\to [0,1]$ with compact support contained in $\Omega_0'\setminus \bar{\Omega }_0$ such that $\chi =1$ on $\overline{\Omega _0''}\setminus \Omega _0'''$. Define a function $\Phi =\tilde{\Phi }-3\lambda \chi $ on $\Omega _0'$, where $\lambda >0$ is chosen small enough so that $\Phi $ is still strictly plurisubharmonic. By choosing $M$ large enough we want to ensure that
\begin{equation}\label{P1}
\Psi (z)>\Phi (z)+2\lambda \quad \textrm{for all } z\in \overline{\Omega _0''}\setminus \Omega _0''' \textrm{ and}
\end{equation}
\begin{equation}\label{P2}
\set{z\in \Omega _0}{\Psi (z)\le \delta } \textrm{ is a compact set for all small enough $\delta >0$}.
\end{equation}
Note that for $z\in Y\cap (\overline{\Omega _0''}\setminus \Omega _0''')$ we have $\Psi (z)=\phi (z)$ and $\Phi (z)=\phi (z)-3\lambda $ (independently of $M$). Thus (\ref{P1}) holds for all $z$ in some open neighborhood $U\subset \CC ^N$ of the set $Y\cap (\overline{\Omega _0''}\setminus \Omega _0''')$ (the bigger the $M$ the bigger the neighborhood). Moreover, the compactness of $(\overline{\Omega _0''}\setminus \Omega _0''')\setminus U$ and the positivity of $\psi $ on this set ensure that (\ref{P1}) will hold for all $M$ big enough. To prove (\ref{P2}) choose an open neighborhood $W\Subset \Omega _0\cup \Omega _1$ of the set $\bar{G}=\set{z\in W'}{\phi (z)\le 0}$ and define $\delta '=\frac{1}{2}\max \set{\phi (z)}{z\in Y\cap \d W}>0$. Note that for $z\in Y\cap \d W$ we have $\Psi (z)=\phi (z)\ge 2\delta '$. Thus $\Psi >\delta '$ on some open neighborhood $U'\subset \CC ^N$ of the set $Y\cap \d W$ (the bigger the $M$ the bigger the neighborhood). Moreover, the compactness  of $\d W\setminus U'$ and the positivity of $\psi $ on this set ensure that for all big enough $M$ we have $\Psi >\delta '$ on $\d W$ and thus (\ref{P2}) holds for all $\delta \le \delta '$.

Now fix an $M$ so that $\Psi $ satisfies (\ref{P1}) and (\ref{P2}). For $\epsilon >0$ let $M_{\epsilon }$ be as in Lemma \ref{sm}. Choose an $\epsilon >0$ such that $\epsilon <\lambda ,\delta '$. Define the function $M_{\epsilon }(\Phi ,\Psi )$ on $\Omega _0''$ as in (v) of Lemma \ref{sm}. Due to (\ref{P1}) and (ii) of Lemma \ref{sm} this function coincides with the function $\Psi -\epsilon $ on $\Omega _0''\setminus \Omega _0'''$ so we can extend it to $\Omega _1\setminus \Omega _0''$ with $\Psi -\epsilon $.

Set $\tilde{\Omega }=\set{z\in \Omega _0\cup \Omega _1}{M_{\epsilon }(\Phi (z),\Psi (z))<0}$. We clearly have $\tilde{\Omega }\subset \set{z\in \Omega _0\cup \Omega _1}{\Psi (z)<\epsilon }$ due to (iii) of Lemma \ref{sm}. Thus (\ref{P2}) ensures that $\tilde{\Omega }\Subset \Omega _0\cup \Omega _1$. Notice that for $z\in Y\cap \Omega _0''$ we have $\Psi (z)=\phi (z)$ and $\Phi (z)\le \phi (z)$, thus $M_{\epsilon }(\Phi (z),\Psi (z))\le \phi (z)$ due to (iv) of Lemma \ref{sm}. Moreover, for $z\in Y\cap (\Omega _1\setminus \Omega _0''')$ we have $\Psi (z)-\epsilon =\phi (z)-\epsilon <\phi (z)$. This implies that $G\subset \tilde{\Omega }$. Lastly, by (iii) of Lemma \ref{sm} we have $M_{\epsilon }(\Phi (z),\Psi (z))\ge \Phi (z)=\tilde{\Phi }(z)$ for $z\in \Omega _0$ and since by construction $r$ maps $\set{z\in \Omega _0}{\tilde{\Phi }(z)<0}$ into $G$ we get $r(\tilde{\Omega }\cap \Omega _0)\subset G$. If the boundary $\d \tilde{\Omega }$ is smooth (i.e., the extended function $M_{\epsilon }(\Phi, \Psi)$ has no critical points on $\d \tilde{\Omega }$) then $\tilde{\Omega }$ is a Stein domain with strongly pseudoconvex boundary of class $\mathcal{C}^k$ and we may take $\Omega =\tilde{\Omega }$. In general this will be true so we need to slightly perturb $\tilde{\Omega }$ to satisfy this last desired property while making sure that the other properties remain valid.

Set $V=\Omega _0''\setminus \supp \chi $. Notice that by construction the function $M_{\epsilon }(\Phi ,\Psi )$ has no critical points on the set $\d \tilde{\Omega }\cap Y\cap \bar{V}$. Hence the same holds true on some compact neighborhood $K\subset \Omega _0''$ of this set. By Morse lemma (see \cite[Chapter VII, Lemma 8.5]{Lau}) we can approximate the extended function $M_{\epsilon }(\Phi ,\Psi )$ arbitrarily closely by another strictly plurisubharmonic function $\rho $ of class $\mathcal{C}^k$ which has only isolated critical points, such that $\rho =M_{\epsilon }(\Phi ,\Psi )$ on $K$. Lastly, to get rid of the critical points on the boundary, we slightly dent the domain $\set{z\in \Omega _0\cup \Omega _1}{\rho (z)<0}$ inwards at each critical point on the boundary (see for example the domain $D_-$ in Main lemma 2.7 and Fig. 1 in \cite{Hen&Lei} on pages 79--80) to get the domain $\Omega $. If the approximation was close enough and the dents were small enough then the domain $\Omega $ still has all the desired properties.
\endproof

For a domain $D$ in a complex manifold $X$ we denote by $\mathcal{C}^l_{0,1}(\bar{D})$ the set of all $(0,1)$-forms whose coefficients are in $\mathcal{C}^l(\bar{D})$. Given an open subset $U\subset D$ we denote by $\mathcal{C}_{0,1}(\bar{D};U)$ the set of all $(0,1)$-forms whose coefficients are in $\mathcal{C}(\bar{D})$ with support in $\bar{U}$. The sets $D$ and $U$ may have some common boundary.

The following proposition gives a bounded linear solution operator for the inhomogeneous $\dbar $-equation for $\dbar $-closed forms from $\mathcal{C}_{0,1}(\bar{D};U)$. The classical Stein version of this result is due to Lieb, Range and Siu (see \cite[p.\ 56, Theorem 2.5.3]{For3} and references therein).

\begin{proposition}\label{bso}
Let $D$ be a relatively compact 1-convex domain with strongly pseudoconvex boundary of class $\mathcal{C}^k$ ($k\ge 2$) in a complex manifold $X$, and let $U\subset D$ be a Stein domain whose closure $\bar{U}$ does not intersect the exceptional set of $D$. Then there exists a linear operator $T:\mathcal{C}_{0,1}(\bar{D};U)\to \mathcal{C}(\bar{D})$ with the following properties:
\begin{enumerate}
\item if $f\in \mathcal{C}_{0,1}(\bar{D};U)\cap \mathcal{C}_{0,1}^1(D)$ and $\dbar f=0$ then $\dbar (Tf)=f$,
\item if $f\in \mathcal{C}_{0,1}(\bar{D};U)\cap \mathcal{C}_{0,1}^l(D)$ for some $l\in \{0,\ldots ,k\}$ then
$$\|Tf\|_{\mathcal{C}^{j,1/2}(\bar{D})}\le M\|f\|_{\mathcal{C}_{0,1}^j(\bar{D})}, \quad j\in \{0,\ldots ,l\}.$$
\end{enumerate}
The constant $M$ depends only on $D$, $U$and $j$.
\end{proposition}

\proof
The idea is to take the Remmert reduction $R:D'\to Y$ of a slightly larger 1-convex domain $D'\supset \bar{D}$, embed the Stein space $Y$ into $\CC ^N$, and then use Proposition \ref{esd} to find a Stein domain $\Omega \subset \CC ^N$ with strongly pseudoconvex boundary such that $\Omega \cap Y=R(D)$, which will allow us to extend forms from $\mathcal{C}_{0,1}(\overline{R(D)};R(U))$ to $\bar{\Omega }$. This will reduce our problem to the Stein case for which the result is already known.

By definition there exist an open neighborhood $W\subset X$ of the boundary $\d D$ and a strictly plurisubharmonic function $\varphi :W\to \RR $ of class $\mathcal{C}^k$ such that $D\cap W=\set{z\in W}{\varphi (z)<0}$ and $\varphi $ has no critical points on $\d D$. Thus, for $\epsilon >0$ small enough, the domain $D':=\set{z\in W}{\varphi (z)<\epsilon }\cup D\subset X$ is also a 1-convex domain whose exceptional set is the same as that of $D$.

Let $R:D'\to Y$ be the Remmert reduction. Denote by $E$ the exceptional set of $D$ and set $G:=R(D)$, $V:=R(U)$, $F:=R(E)$. Then $G$ is a Stein domain with strongly pseudoconvex boundary of class $\mathcal{C}^k$, $V$ is a Stein domain, and $Y$ has only isolated singularities contained in $F$. Each form $\alpha \in \mathcal{C}_{0,1}^j(\bar{D};U)$ can be naturally identified with a form in $\mathcal{C}_{0,1}^j(\bar{G};V)$. Indeed, since $R:D'\setminus E\to Y\setminus F$ is a biholomorphism, the push-forward $R_{\star }\alpha $ is well defined on $\bar{G}\setminus F$. Moreover, since $\alpha $ has support in $\bar{U}$, $R_{\star }\alpha $ has support in $\bar{V}$ and can therefore be extended with $0$ across $F$ (since $F\cap \bar{V}=\emptyset $). Thus $R_{\star }:\mathcal{C}_{0,1}(\bar{D};U)\to \mathcal{C}_{0,1}(\bar{G};V)$ is a well defined linear operator which commutes with $\dbar $. Moreover, it preserves norms and degree of smoothness of the forms.

By classical results (see \cite{Bis}, \cite{Nar}) there exists a proper holomorphic embedding $i:Y\to \CC ^N$ (for some $N$). We shall identify $Y$ with its image $i(Y)$. According to the tubular neighborhood theorem due to Docquier and Grauert (see \cite{Doc&Gra} or \cite[p.\ 67, Theorem 3.3.3]{For3}) there exists an open (tubular) neighborhood $\Omega _0'\subset \CC ^N$ of $\bar{V}$ and a holomorphic retraction $r:\Omega _0'\to Y\cap \Omega _0'$. Choose an open neighborhood $\Omega _0\Subset \Omega _0'$ of $\bar{V}$ and an open neighborhood $\Omega _1\subset \CC ^N$ of $\bar{G}\setminus \Omega _0$ such that $r(\Omega _1\cap \Omega_0)\cap \bar{V}=\emptyset $.

By Proposition \ref{esd} there exists a Stein domain $\Omega \Subset \CC ^N$ with strongly pseudoconvex boundary of class $\mathcal{C}^k$ such that $r(\Omega \cap \Omega _0)\subset G\subset \Omega \Subset \Omega _0\cup \Omega _1$.

Every $(0,1)$-form $\alpha \in \mathcal{C}_{0,1}^j(\bar{G};V)$ can be extended to a $(0,1)$-form $P(\alpha )\in \mathcal{C}_{0,1}^j(\bar{\Omega };V')$, where $V'=r^{-1}(V)\cap \Omega $, by setting $P(\alpha )=r^{\star }\alpha $ (the pull-back) on $\bar{\Omega }\cap \Omega _0$ and $P(\alpha )=0$ on $\bar{\Omega }\cap \Omega _1$. Indeed, by the choice of $\Omega _1$ the definitions coincide on the intersection $\bar{\Omega }\cap \Omega_0\cap \Omega _1$ and since $r$ is holomorphic $P(\alpha )$ is a $(0,1)$-form. Notice that $P:\mathcal{C}_{0,1}(\bar{G};V)\to \mathcal{C}_{0,1}(\bar{\Omega };V')$ is a linear operator that commutes with $\dbar $. Moreover, it preserves norms and degree of smoothness of the forms. Now we are in a position to use a result due to Lieb, Range and Siu (see \cite[p.\ 56, Theorem 2.5.3]{For3}) which gives us a bounded linear solution operator to the $\dbar $-equation $\tilde{T}:\mathcal{C}_{0,1}(\bar{\Omega })\to \mathcal{C}(\bar{\Omega })$. Let $S:\mathcal{C}(\bar{\Omega })\to \mathcal{C}(\bar{G})$ be the restriction operator $S(\alpha )=\alpha |_{\bar{G}}$. 
The operator 
$$T=R^{\star }\circ S\circ \tilde{T}\circ P\circ R_{\star }:\mathcal{C}_{0,1}(\bar{D};U)\to \mathcal{C}(\bar{D})$$
then enjoys the desired properties.
\endproof

\begin{remark}
It would be of some interest to know if a solution operator with the analogous properties could be found for $(0,1)$-forms whose coefficients vanish only in some small fixed neighborhood of the exceptional set of $D$, or even for $(0,1)$-forms whose coefficients vanish only along the exceptional set of $D$ to some high enough order. To prove the proposition for the forms of the first kind using the same method we would need to have a tubular neighborhood of the set $\bar{G}\setminus W$, where $W$ is a neighborhood of the set $F$. However, according to \cite[Theorem 7.1]{Ros} this would actually mean that $G$ has no singularities. So this would only work for 1-convex domains $D$ whose Remmert reduction is a Stein manifold.
\end{remark}

\section{A generalized Cartan splitting lemma}
\label{sec:Cartan}

In this section we prove a version of the generalized Cartan lemma for 1-convex Cartan pairs with estimates up to the boundary (see Theorem \ref{st}). The definition of a 1-convex Cartan pair is the following (compare with Definitions 2.1 and 2.6 in \cite{Hen&Lei}; for the definition of a Cartan pair in the standard Stein case see \cite[p.\ 209, Definition 5.7.1]{For3}). 

\begin{definition}\label{cp}
A pair $(D_0,D_1)$ of domains in a complex manifold $X$ is a \emph{1-convex Cartan pair} of class $\mathcal{C}^k$ if it satisfies the following conditions:
\begin{enumerate}
\item $D_0$, $D_1$ and $D_0\cup D_1$ are relatively compact 1-convex domains with strongly pseudoconvex boundary of class $\mathcal{C}^k$, and $D_0\cap D_1$ is a relatively compact Stein domain with strongly pseudoconvex boundary of class $\mathcal{C}^k$,
\item $\overline{D_0\setminus D_1}\cap \overline{D_1\setminus D_0}=\emptyset $,
\item the exceptional set $E$ of $D_0\cup D_1$ does not intersect $\overline{D_0\cap D_1}$.
\end{enumerate}
We say that $D_1$ is a \emph{convex bump} on $D_0$ if, in addition to the above, there exists a biholomorphic map from an open neighborhood of $\bar{D}_1$ onto an open set in $\CC ^n$ ($n=\dim _{\CC }X$) that maps $D_1$ and $D_0\cap D_1$ onto strongly convex domains.
\end{definition}

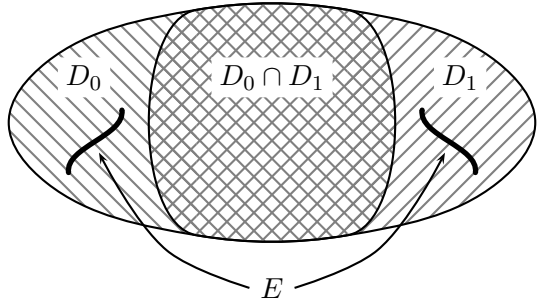
\begin{figure}[h]

\begin{pspicture}(-3.4,-2.3)(3.4,1.7)

\psecurve[fillstyle=vlines,hatchcolor=gray](3,0.8)(1,1.5)(-1,1.5)(-3,0.8)(-3.5,0)(-3,-0.8)(-1,-1.5)(1,-1.5)(3,-0.8)
\psecurve[fillstyle=hlines,hatchcolor=gray](-3,-0.8)(-1,-1.5)(1,-1.5)(3,-0.8)(3.5,0)(3,0.8)(1,1.5)(-1,1.5)(-3,0.8)
\psecurve[fillstyle=crosshatch,hatchcolor=gray](0,-1.5)(-1,1.5)(-1,-1.5)(0,1.5)
\psecurve[fillstyle=crosshatch,hatchcolor=gray](0,-1.5)(1,1.5)(1,-1.5)(0,1.5)

\psecurve[linewidth=2pt,arrows=cc-cc](-2.3,-1)(-2.7,-0.7)(-2,0.2)(-2.4,0.5)
\psecurve[linewidth=2pt,arrows=cc-cc](2.3,-1)(2.7,-0.7)(2,0.2)(2.4,0.5)

\rput*(-2.5,0.6){$D_0$}
\rput*(2.5,0.6){$D_1$}
\rput*(0,0.6){$D_0\cap D_1$}
\rput(0,-2.2){$E$}
\pscurve[arrows=->](-0.3,-2.2)(-1.4,-1.7)(-2.3,-0.4)
\pscurve[arrows=->](0.3,-2.2)(1.4,-1.7)(2.3,-0.4)

\end{pspicture}

\caption{A 1-convex Cartan pair $(D_0,D_1)$.}
\label{figcp}
\end{figure}

Note that the set $D_1$ in the definition of a convex bump is of course Stein.

The following theorem is a 1-convex version of the generalized Cartan splitting lemma found in \cite[Theorem 3.2]{Dri&For} (see also \cite[Proposition 5.8.1]{For3}).

\begin{theorem}\label{st}
Let $(D_0,D_1)$ be a 1-convex Cartan pair of class $\mathcal{C}^l$ $(l\ge 2)$ in a complex manifold $X$. Set $D_{0,1}=D_0\cap D_1$ and $D=D_0\cup D_1$. Given a bounded open convex set $P\in \CC ^N$ containing the origin and a number $r\in (0,1)$ there exists a number $\epsilon >0$ satisfying the following. For every map $\gamma :\bar{D}_{0,1}\times P\to \CC ^N$ of class $\AA^k$ $(k\in\{0,1\ldots, l\})$ of the form
$$\gamma (x,t)=t+c(x,t), \quad x\in \bar{D}_{0,1}, t\in P,$$
with $\|c\|_{\mathcal{C}^k(\bar{D}_{0,1}\times P)}<\epsilon $, there exist maps $\alpha :\bar{D}_0\times rP\to \CC ^N$ and $\beta :\bar{D}_1\times rP\to \CC ^N$ of class $\AA ^k$ and of the form
$$\alpha (x,t)=t+a(x,t), \quad x\in \bar{D}_0,\ t\in rP,$$
$$\beta (x,t)=t+b(x,t), \quad x\in \bar{D}_1,\ t\in rP,$$
depending smoothly on $\gamma $, satisfying the condition
$$\gamma (x,\alpha (x,t))=\beta (x,t), \quad x\in \bar{D}_{0,1},\ t\in rP$$
and the estimates
$$\|a\|_{\mathcal{C}^k(\bar{D}_0\times rP)}\le M\|c\|_{\mathcal{C}^k(\bar{D}_{0,1}\times P)},$$
$$\|b\|_{\mathcal{C}^k(\bar{D}_1\times rP)}\le M\|c\|_{\mathcal{C}^k(\bar{D}_{0,1}\times P)},$$
where the constant $M$ depends only on $D$. If $c$ vanishes to order $m\in \ZZ _+$ along $t=0$ then so do $a$ and $b$. Moreover, if $X'$ is a closed analytic subset of $X$ that does not intersect $\bar{D}_{0,1}$ we can insure that $a$ and $b$ vanish to any given finite order along $(X'\cap \bar{D}_0)\times rP$ and $(X'\cap \bar{D}_1)\times rP$, respectively.
\end{theorem}

Denote by $A_r$, $B_r$ and $C_r$ the Banach spaces of maps $\bar{D}_0\times rP\to \CC ^N$, $\bar{D}_1\times rP\to \CC ^N$ and $\bar{D}_{0,1}\times rP\to \CC ^n$, respectively, which are of class $\AA ^k$ and have finite $\mathcal{C}^k$ norm.

We first prove the following lemma (cf.\ \cite[p.\ 212, Lemma 5.8.2]{For3}).

\begin{lemma}\label{sl}
Assume the hypotheses of Theorem \ref{st}. There exist bounded linear operators $\AA :C_r\to A_r$ and $\mathcal{B}:C_r\to B_r$ satisfying
$$\phi =\AA \phi -\mathcal{B}\phi , \quad \phi \in C_r. $$
If $\phi $ vanishes to order $m\in \ZZ _+$ along $t=0$ then so do $\AA \phi $ and $\mathcal{B}\phi $. If moreover $X'\subset X$ is a closed analytic subset that does not intersect $\bar{D}_{0,1}$, then we can insure that $\AA \phi $ and $\mathcal{B}\phi $ agree with $Id$ to any given finite order along $(X'\cap \bar{D}_0)\times rP$ and $(X'\cap \bar{D}_1)\times rP$, respectively.
\end{lemma}

\proof
The proof is essentially the same as in \cite[p.\ 212]{For3}. Because of the condition (ii) in Definition \ref{cp} there is a smooth function $\chi :X\to [0,1]$ such that $\chi =0$ in a neighborhood of $\overline{D_0\setminus D_1}$ and $\chi =1$ in a neighborhood of $\overline{D_1\setminus D_0}$. Thus for any $\phi \in C_r$ the product $\chi (x)\,\phi (x,t)$ extends to a $\mathcal{C}^k$ function on $\bar{D}_0\times rP$ which vanishes on $\overline{D_0\setminus D_1}\times rP$, and the product $(\chi (x)-1)\, \phi (x,t)$ extends to a $\mathcal{C}^k$ function on $\bar{D}_1\times rP$ which vanishes on $\overline{D_1\setminus D_0}\times rP$. Moreover, we have
\begin{eqnarray*}
\dbar _x(\chi \phi (\cdot ,t)) & = & \phi (\cdot ,t)\dbar \chi \qquad \textrm{on } D_0, \\
\dbar _x((\chi -1)\phi (\cdot ,t)) & = & \phi (\cdot ,t)\dbar \chi \qquad \textrm{on } D_1.
\end{eqnarray*}
Here $\phi (\cdot ,t)\,\dbar \chi $ extends to a $(0,1)$-form on $\bar{D}$ of class $\mathcal{C}^k$ which depends holomorphically on the parameter $t$.

Choose functions $f_1,\ldots ,f_l\in \OO (X)$ that vanish to order $m$ along the subvariety $X'$ and have no common zeros on $\bar{D}_{0,1}$. Since $D_{0,1}$ is by assumption a Stein domain with strongly pseudoconvex boundary, Cartan's division theorem (see \cite[p.\ 54, Corollary 2.4.4]{For3}) gives holomorphic functions $g_1,\ldots ,g_l$ in some neighborhood $U$ of $\bar{D}_{0,1}$ such that $\sum _{j=1}^lf_jg_j=1$ on $U$. Notice that $g_j\phi (\cdot ,t)\dbar \chi $ is a family of $\dbar $-closed $(0,1)$-forms on $\bar{D}$ with support in $\bar{D}_{0,1}$ that depend holomorphically on the parameter $t$. Moreover, by the condition (iii) in the Definition \ref{cp} $\bar{D}_{0,1}$ does not intersect the exceptional set of $D$. According to Proposition \ref{bso} there exists a bounded linear solution operator $T:\mathcal{C}_{0,1}(\bar{D};D_{0,1})\to \mathcal{C}(\bar{D})$ to the $\dbar $-equation. For $\phi \in C_r$ we set
for every  $t\in rP$:
\begin{eqnarray*}
(\mathcal{A}\phi )(x,t) & := & \chi (x)\, \phi (x,t)-\sum _{j=1}^lf_j(x)\, T(g_j\phi (\cdot ,t)\, \dbar \chi )(x), \quad \hspace{1cm} x\in \bar{D}_0,\\
(\mathcal{B}\phi )(x,t) & := & (\chi (x)-1)\, \phi (x,t)-\sum _{j=1}^lf_j(x)\, T(g_j\phi (\cdot ,t)\, \dbar \chi )(x), \quad x\in \bar{D}_1.
\end{eqnarray*}
Clearly $\phi =\AA \phi -\mathcal{B}\phi $ on $\bar{D}_{0,1}\times rP$. Since $T$ is the solution operator to the $\dbar $-equation (with respect to variable $x$) and it commutes with derivations with respect to variables $t\in \CC ^N$, we have $\dbar _x(\AA \phi )=0$, $\dbar _x(\mathcal{B}\phi )=0$ and $\dbar _t(\AA \phi )=0$, $\dbar _t(\mathcal{B}\phi )=0$ in the interior of the respective domains. Also, the choice of the functions $f_1,\ldots ,f_l$ ensures the vanishing of $\AA \phi $ and $\mathcal{B}\phi $ to order $m$ along $(X'\cap \bar{D}_0)\times rP$ and $(X'\cap \bar{D}_1)\times rP$, respectively. The boundedness of $\AA $ and $\mathcal{B}$ follows from  boundedness of $T$.
\endproof

\proof[Proof of Theorem \ref{st}]
Fix a number $R\in (r,1)$ and set $\gamma _0(x,t)=t$. Let $\AA$ and $\mathcal{B}$ be as in Lemma \ref{sl}. For $\gamma \in C_1$ close to $\gamma _0$ and $\phi \in C_r$ close to $0$ we define a map
$\Psi (\gamma,\phi)$ with values in $\mathbb{C}^N$ by setting
$$\Psi (\gamma ,\phi )(x,t)=\gamma (x,t+\AA \phi (x,t))-(t+\mathcal{B}\phi (x,t)), \quad x\in \bar{D}_{0,1}, t\in rP.$$
Then $(\gamma ,\phi )\mapsto \Psi (\gamma ,\phi )$ is a smooth map from an open neighborhood of $(\gamma _0,0)$ in the Banach space $C_1\times C_r$ to the Banach space $C_r$. Indeed, $\Psi $ is obviously continuous, it is linear in the variable $\gamma $, and its partial derivative with respect to $\phi$ equals
$$(\d _{\phi }\Psi (\gamma ,\phi )\tau )(x,t)=\d _t\gamma (x,t+\AA \phi (x,t))\cdot \AA \tau (x,t)-\mathcal{B}\tau (x,t).$$
This is again linear in $\gamma $ and continuous due to the following Cauchy estimates in the variable $t$,
$$\sup \set{|\d _x^{\mu }\d _t^{\nu +1}\gamma (x,t)|}{|\mu |+|\nu |\le k, x\in \bar{D}_{0,1}, t\in RP}\le C\|\gamma \|_{\mathcal{C}^k(\bar{D}_{0,1}\times P)},$$
which take care of the $\d _t\gamma $ part. Similar argument applies to higher order derivatives of $\Psi$.

According to Lemma \ref{sl} we have
$$\Psi (\gamma _0,\phi )=\AA \phi -\mathcal{B}\phi =\phi ,$$
hence $\d _\phi \Psi (\gamma _0,\phi )$ is the identity map.\ By the implicit function theorem there exists a smooth map $\gamma \mapsto \Phi (\gamma )$ from an open neighborhood of $\gamma _0$ in $C_R$ to an open neighborhood of $0$ in $C_r$ satisfying $\Phi (\gamma _0)=0$ and $\Psi (\gamma ,\Phi (\gamma ))=0$. The maps
$$\alpha _{\gamma }(x,t)=t+\AA \circ \Phi (\gamma )(x,t), \quad \beta _{\gamma }(x,t)=t+\mathcal{B}\circ \Phi (\gamma )(x,t)$$
satisfy all the conclusions in the theorem.
\endproof

\section{Gluing method}
\label{sec:gluing}

The following proposition is what we have in mind when we talk about gluing sprays of sections. It is one of the most important tools in this theory. The standard Stein case was obtained in \cite[Proposition 4.3]{Dri&For} (see also \cite[Proposition 2.4]{Dri&For1} and \cite[p.\ 216, Proposition 5.9.2]{For3}).

\begin{proposition}\label{gs}
Let $(D_0,D_1)$ be a 1-convex Cartan pair of class $\mathcal{C}^l$ $(l\ge 2)$ in a complex manifold $X$. Set $D_{0,1}=D_0\cap D_1$ and $D=D_0\cup D_1$. Let $\pi :Z\to \bar{D}$ be either a fiber bundle of class $\AA ^k$ $(k\le l)$ or the restriction to $\bar D$ of a holomorphic submersion $\tilde{Z}\to X$. Given a spray of $\pi $-sections $f:\bar{D}_0\times P_0\to Z$ of class $\AA ^k$ which is dominating on $\bar{D}_{0,1}$, there exists an open set $0\in P\subset P_0$ such that the following hold:
\begin{enumerate}
\item for every spray of $\pi $-sections $g:\bar{D}_1\times P_0\to Z$ of class $\AA ^k$ which is sufficiently close to $f$ in $\mathcal{C}^k(\bar{D}_{0,1}\times P_0)$ there exists a spray of $\pi $-sections $\tilde{f}:\bar{D}\times P\to Z$ of class $\AA ^k$ which is close to $f$ in $\mathcal{C}^k(\bar{D}_0\times P)$ (with respect to the $\mathcal{C}^k(\bar{D}_{0,1}\times P_0)$ distance between $f$ and $g$), whose core $\tilde{f}_0$ is homotopic to $f_0$ on $\bar{D}_0$ and to $g_0$ on $\bar{D}_1$,
\item if $f$ and $g$ agree to order $m\in \ZZ _+$ along $\bar{D}_{0,1}\times \{0\}$ then we can ensure that $\tilde{f}$ agrees to order $m$ with $f$ along $\bar{D}_0\times \{0\}$ and with $g$ along $\bar{D}_1\times \{0\}$, and
\item if $\sigma $ is a zero set of finitely many $\AA ^k(D_0)$ functions and $\sigma \cap \bar{D}_{0,1}=\emptyset $ then we can ensure that $\tilde{f}_0$ agrees with $f_0$ to order $m$ along $\sigma $.
\end{enumerate}
\end{proposition}

We begin with a lemma which gives us a transition map between a pair of nearby sprays of sections. (Compare with \cite[Proposition 4.4]{Dri&For} and \cite[p.\ 216, Proposition 5.9.3]{For3}.)

\begin{lemma}\label{tm}
Assume the hypotheses of Proposition \ref{gs}. Let $\epsilon >0$. There exists an open (convex) set $P_1\subset P_0$ containing the origin such that if $f$ and $g$ are sufficiently close in $\mathcal{C}^k(\bar{D}_{0,1}\times P_0)$ there exists a map $\gamma :\bar{D}_{0,1}\times P_1\to \CC ^N$ of class $\AA ^k$ satisfying
$$\gamma (x,t)=t+c(x,t),\quad \|c\|_{\mathcal{C}^k(\bar{D}_{0,1}\times P_1)}<\epsilon ,$$
$$f(x,t)=g(x,\gamma (x,t)), \quad x\in \bar{D}_{0,1},\ t\in P_1.$$
If $f$ and $g$ agree to order $m\in \ZZ _+$ along $\bar{D}_{0,1}\times \{0\}$ then $\gamma $ can be chosen of the form $\gamma (x,t)=t+\sum _{|J|=m}c_J(x,t)t^J$ with $c_J\in \AA ^k(D_{0,1}\times P_1)^N$. 
\end{lemma}

\proof
Denote by $E$ the subbundle of $\bar{D}_{0,1}\times \CC ^N$ with fibers
$$E_x=\ker (\d _t|_{t=0}f(x,t):\CC ^N\to VT_{f(x,0)}Z),\quad x\in \bar{D}_{0,1}.$$
Notice that $E$ is a fiber bundle of class $\AA ^k$. By Theorem B for such fiber bundles, due to Heunemann  \cite{Heu}, there exists a subbundle $E'$ of $\bar{D}_{0,1}\times \CC ^N$ of class $\AA ^k$ such that $\bar{D}_{0,1}\times \CC ^N=E\oplus E'$. For each $x\in \bar{D}_{0,1}$ and $t\in \CC ^N$ we write $t=t_x\oplus t'_x\in E\oplus E'$. The map
$$\d _t|_{t=0}f(x,t):E'\to VTZ|_{f_0(\bar{D}_{0,1})}$$
is an isomorphism. By the implicit function theorem there is an open (convex) set $P_1\subset P_0$ containing the origin such that for every spray of $\pi $-sections $g:\bar{D}_1\times P_0\to Z$ of class $\AA ^k$ which is sufficiently close to $f$ in $\mathcal{C}^k(\bar{D}_{0,1}\times P_0)$ there is a unique map
$$\tilde{\gamma }(x,t)=\tilde{\gamma }(x,t_x\oplus t'_x)=t_x\oplus (t'_x+\tilde{c}(x,t))\in E_x\oplus E'_x\cong \CC ^N$$
of class $\AA ^k(\bar{D}_{0,1}\times P_1)$ satisfying $f(x,\tilde{\gamma }(x,t))=g(x,t)$, where $\|\tilde{c}\|_{\mathcal{C}^k(\bar{D}_{0,1}\times P_1)}$ is controlled by the $\mathcal{C}^k(\bar{D}_{0,1}\times P_0)$ distance between $f$ and $g$. Shrinking $P_1$ if necessary, the map $\tilde{\gamma }$ has a fiberwise inverse $\gamma (x,t)=t_x\oplus (t'_x+\hat{c}(x,t))=t+c(x,t)$, which satisfies all the conclusions in the lemma.
\endproof

\proof[Proof of Proposition \ref{gs}]
Let the open (convex) set $P_1\subset P_0$ and the map $\gamma (x,t)=t+c(x,t)$ be as in Lemma \ref{tm}. Let $P=rP_1$ with $r\in (0,1)$. If $g$ is sufficiently close to $f$ in $\mathcal{C}^k(\bar{D}_{0,1}\times P_0)$ then $c$ is sufficiently close to $0$ in $\mathcal{C}^k(\bar{D}_{0,1}\times P_1)$ and by Theorem \ref{st} there exist maps $\alpha :\bar{D}_0\times P\to \CC ^N$ and $\beta :\bar{D}_1\times P\to \CC ^N$ of class $\AA ^k$ satisfying
$$\gamma (x,\alpha (x,t))=\beta (x,t), \quad x\in \bar{D}_{0,1}, t\in P.$$
Combining this with the equality
$$f(x,t)=g(x,\gamma (x,t)), \quad x\in \bar{D}_{0,1}, t\in P_1$$
from Lemma \ref{tm} we get
$$f(x,\alpha (x,t))=g(x,\gamma (x,\alpha (x,t)))=g(x,\beta (x,t)), \quad x\in \bar{D}_{0,1}, t\in P.$$
Therefore the maps $f(x,\alpha (x,t))$ and $g(x,\beta (x,t))$ amalgamate into a spray of $\pi $-sections $\tilde{f}:\bar{D}\times P\to Z$ which has the desired properties.
\endproof

\section{Existence of sprays}
\label{sec:existence}

In this section we prove Theorem \ref{eost}.

The most common way to obtain sprays, and sprays of sections, is the following. Suppose that $Z|_D$ admits finitely many vertical holomorphic vector fields $V_j$, $j=1,\ldots ,N$ (i.e., holomorphic sections of the bundle $VTZ|_D$), whose flows $\phi _t^j$ exist (and are holomorphic) for all sufficiently small times $t\in P\subset \CC $. Then the map $s:Z|_D\times P^N\to Z$ defined by
$$s(z;t_1,\ldots ,t_N)=\phi _{t_N}^N\circ \ldots \circ \phi _{t_1}^1(z)$$
is a $\pi $-spray on $Z|_D$. We have $\d _{t_j}s(z;0,\ldots ,0)=V_j(z)$. If the vector fields $V_j$ span $VTZ$ at each point $z\in K\subset Z|_D$ then the spray is dominating on $K$. Similarly, for every holomorphic $\pi $-section $f_0:D\to Z$ the map $f:D\times P^N\to Z$, defined by
$$f(x;t_1,\ldots ,t_N)=s(f_0(x);t_1,\ldots ,t_N)=\phi _{t_N}^N\circ \ldots \circ \phi _{t_1}^1\circ f_0(x),$$
is a holomorphic spray of $\pi $-sections with core $f_0$ which is dominating on $f_0^{-1}(K)$. Moreover, if all the flows $\phi _t^j$ are stationary on $L\subset Z|_D$ then $f$ fixes $f_0^{-1}(L)$. If the vector fields $V_j$ are complete (in the sense that their flows exist for all times $t\in \CC $) then the sprays defined above are global. For the construction of $f$ we only need vector fields $V_j$ to be defined on some neighborhood of the graph $f_0(D)\subset Z|_D$. In the case where $D$ is a Stein domain such vector fields can always be constructed using Cartan's theorem $A$ applied on a Stein neighborhood of the set $f_0(D)$. Hence we have the following proposition. (For a stronger version and a detailed proof see \cite[p.\ 220, Lemma 5.10.4]{For3}; for the case with boundary see \cite[Corollary 4.2]{Dri&For1} and \cite[Lemma 4.2]{Dri&For}.)

\begin{proposition}\label{eosp}
Let $D$ be a relatively compact domain in a Stein manifold $S$ and $\pi :Z\to S$ a holomorphic submersion. Given a holomorphic $\pi$-section $f_0:S\to Z$ there exists a holomorphic spray of $\pi$-sections $f:D\times P\to Z$ with core $f_0|_{D}$ which is dominating on $D$.
\end{proposition}

Assume now that $X$ is a complex manifold and $D\Subset X$ is a 1-convex domain with $\mathcal{C}^l$ ($l\ge 2$) boundary and exceptional set $E$. Suppose $\pi :Z\to \bar{D}$ is either a fiber bundle of class $\AA ^k$ ($k\le l$) or the restriction of a holomorphic submersion $\tilde{Z}\to X$. In either case each point $z\in Z$ admits an open neighborhood $\Omega \subset Z$ isomorphic to $U\times V$ where $U$ is a relatively open subset of $\bar{D}$ and $V$ is an open subset of $\CC ^m$, such that in coordinates $z=(x,y)\in U\times V$ the map $\pi $ is just the projection $(x,y)\to x$. Such $\Omega $ will be called a \emph{special coordinate chart} on $Z$.

Given a $\pi $-section $f_0:\bar{D}\to Z$ of class $\AA ^k$ we would like to find a spray of $\pi $-sections $f:\bar{D}\times P\to Z$ of class $\AA ^k$ with core $f_0$ which is dominating on $\bar{D}\setminus E$ and fixes $E$. The main problem is to find such a spray in some neighborhood of the exceptional set $E$.

An attempt in this direction was made in \cite[\S 4]{Pre} (see also \cite[Corollary 2.6]{Pre&Sla}) by following the method described above. In this case the set $f_0(D)$ of course does not have a Stein neighborhood (in most cases it does not even have a 1-convex neighborhood). Hence we have to consider the set $D\setminus g^{-1}(0)$ instead of $D$, where $g:D\to \CC $ is a holomorphic function which vanishes on $E$. The set $D\setminus g^{-1}(0)$ is then Stein and we have a Stein neighborhood $\Omega \subset Z$ for $f_0(D\setminus g^{-1}(0))$. However, given a vector field $V$ on $\Omega $ as above, its flow may not exist (stay in $\Omega $) for all small times since $\Omega $ may be very thin (in the fiber direction) close to $f(g^{-1}(0))$. One of the main contributions of \cite{Pre} is the construction of neighborhoods $\Omega $ that are \emph{conic} along $f_0(g^{-1}(0))$, meaning that the width of $\Omega$ in the fiber direction decreases at most polynomially with the distance from $f_0(g^{-1}(0))$. This is still not enough for the flow to be defined for all small times. Additionally, we would need the vector field $V$ to grow at most polynomially when approaching the boundary of $\Omega $; this is what is missing in \cite{Pre}. Having such a vector field we would then consider the vector field $G^nV$, where $G$ is a fiberwise constant extension of the function $g$ (which we think of as a function on $f_0(D)$). For $n$ big enough this vector field would then have the flow defined for all small times and this flow could be extended across $f_0(g^{-1}(0))$ by identity. A finite number of such vector fields would generate the vertical tangent bundle $VTZ$ on $\Omega $. Hence, the corresponding spray of sections would be dominating outside $g^{-1}(0)$ and would fix $g^{-1}(0)$. With some additional work we could correct such a spray to be dominating outside $E$ and fix $E$.

At the moment we do not know how to construct vector fields with at most polynomial growth. To avoid this difficulty we have introduced Condition $\mathcal{E}$ (see \S 2) which explicitly guarantees the existence of sprays of sections with desired properties in small neighborhoods of the exceptional set $E$.

The most important sufficient conditions implying Condition $\mathcal{E}$ are given by Proposition \ref{ce} which we now prove.

\proof[Proof of Proposition \ref{ce}]
(i) Let $U_0\subset X$ be an open neighborhood of $E$ and $f_0:U_0\to Z$ a holomorphic $\pi $-section. Choose an open 1-convex neighborhood $U\Subset U_0\cap D$ of $E$. We may consider $V|_{f_0(U)}$ as a holomorphic vector bundle over $U$. 

We claim that there exist finitely many holomorphic sections $W_1,\ldots ,W_N$ of the vector bundle $V|_{f_0(U)}\to U$ which span $V_{f_0(x)}$ for each $x\in U\setminus E$ and vanish on $E$. This can be seen by the following argument. Denote by $\mathcal{F}$ the coherent analytic sheaf of germs of sections of the vector bundle $V|_{f_0(U)}\to U$ and let $R:U\to S$ be the Remmert reduction of $U$. Since $R$ is a proper holomorphic map, the direct image sheaf $R_*\mathcal{F}$ is also a coherent analytic sheaf by Grauert's direct image theorem. Since $S$ is Stein, we can apply Cartan's Theorem A to find a finite number of global sections $\tilde{s}_1,\ldots ,\tilde{s}_N\in R_*\mathcal{F}(U)$ which locally generate $R_*\mathcal{F}$ on $U$ and vanish on the finite set $R(E)$. These sections give rise to the corresponding global sections $s_1,\ldots ,s_N\in \mathcal{F}(U)$. The restriction $R:U\setminus E\to S\setminus R(E)$ of the Remmert reduction  is a biholomorphism. Hence, for each $x\in U\setminus E$ the stalk $\mathcal{F}_x$ is isomorphic to the stalk $(R_*\mathcal{F})_{R(x)}$. This implies that the sections $s_1,\ldots ,s_N$ locally generate $\mathcal{F}$ on $U\setminus E$. These sections correspond to holomorphic sections $W_1,\ldots ,W_N$ of the vector bundle $V|_{f_0(U)}\to U$ which span $V_{f_0(x)}$ for each $x\in U\setminus E$ and vanish on $E$.

For a sufficiently small open set $0\in P\subset \CC ^N$ the map
$$f(x,t)=s\Big(\sum _{j=1}^Nt_jW_j(x)\Big), \quad x\in U,\ t\in P$$
is a holomorphic spray of $\pi $-sections with core $f_0|_U$. We have
$$\d _{t_j}|_{t=0}f(x,t)=Ds(0_{f_0(x)})W_j(x).$$
Since the vector fields $W_1,\ldots ,W_N$ span $V_{f_0(x)}$ at each point $x\in U\setminus E$ and the $\pi $-spray $s$ is dominating on $Z|_{D\setminus E}$ we see that $f$ is dominating on $D\setminus E$. Moreover, $f$ fixes $E$ since $W_1,\ldots ,W_j$ vanish on $E$.

\medskip
(ii) By definition, a complex manifold $Y$ is elliptic if (and only if) it admits a global dominating spray in the sense of Definition \ref{sg} with respect to the constant map $Y\to \{*\}$. It is then clear that the canonical projection $X\times Y\to X$ also admits a global dominating spray. Therefore (ii) is just a special case of (i).

\medskip
(iii) Let $U_0\subset X$ be an open neighborhood of $E$ and $f_0:U_0\to Y$ a holomorphic map.\ Fix an open neighborhood $U\Subset U_0$ of $E$. Since $f_0(\bar{U})\subset Y$ is a compact set, the assumption on $Y$ gives finitely many global holomorphic vector fields $V_1,\ldots ,V_N$ on $Y$ that span the tangent space $T_yY$ at each point $y\in f_0(\bar{U})$. Denote by $\phi _t^j$ the flow of $V_j$. The map $f:U\times P^N\to Y$ defined by
$$f(x;t_1,\ldots ,t_N)=\phi _{t_N}^j\circ \ldots \circ \phi _{t_1}^1(f_0(x))$$
is a holomorphic spray of maps with core $f_0|_U$. We have
$$\d _{t_j}|_{t=0}f(x,t)=V_j(f_0(x)).$$
By the choice of the vector fields this means that $f$ is dominating on $U$; however, it does not necessarily fix $E$. To accommodate for this we choose finitely many holomorphic functions $g_1,\ldots ,g_M\in \OO (U)$ such that their common zero set is exactly $E$. We then replace each term $\phi _{t_j}^j$ in the definition of $f$ by a composition of $M$ terms $\phi _{t_{j,k}g_k(x)}^j$ for $k=1,\ldots ,M$. The new spray of maps is then dominating on $U\setminus E$ and fixes $E$.
\endproof

We now turn to the proof of Theorem \ref{eost}. We shall use the standard technique of `extending' the spray of sections step by step to bigger and bigger domains in the base manifold until we reach $\bar{D}$. The initial spray of sections is provided by Condition $\mathcal{E}$. The main step in the proof will be to show how to `extend' a spray of sections across a 1-convex Cartan pair. By `extending' we actually mean to approximate a spray of sections, defined over a smaller set in the base manifold, by a spray of sections defined over a bigger set in the base manifold. We follow the proof in \cite{Dri&For1} which pertains to the Stein case. In each step we will need to ensure that the new spray of sections has core $f_0$, it fixes $E$, and (most importantly) is still dominating outside $E$. The parameter set $0\in P\subset \CC^N$ of the sprays of sections will be allowed to shrink around $0$ at every step.

\begin{proposition}\label{eacp}
Let $(D_0,D_1)$ be a 1-convex Cartan pair of class $\mathcal{C}^l$ $(l\ge 2)$ in a complex manifold $X$ such that $\bar{D}_1$ is contained in some coordinate chart on $X$. Set $D_{0,1}=D_0\cap D_1$ and $D=D_0\cup D_1$ and denote by $E$ the exceptional set of $D_0$. Let $\pi:Z\to \bar{D}$ be either a fiber bundle of class $\AA ^k$ ($k\le l$) or the restriction to $\bar D$ of a holomorphic submersion $\tilde{Z}\to X$.
Let $f_0:\bar{D}\to Z$ be a $\pi $-section of class $\AA ^k$ and $f:\bar{D}_0\times P_0\to Z$ ($0\in P_0\subset \CC ^N$) a spray of $\pi $-sections of class $\AA ^k$ with core $f_0|_{\bar{D}_0}$ which is dominating on $\bar{D}_0\setminus E$ and fixes $E$. Suppose that $f_0(\bar{D}_1)$ is contained in a special coordinate chart on $Z$. If $N\ge \dim _{\CC }Z$ then there exist an open set $0\in P\subset P_0$ and a spray of $\pi $-sections $F:\bar{D}\times P\to Z$ of class $\AA ^k$, with core $f_0$, which is dominating on $\bar{D}\setminus E$ and fixes $E$. Moreover, $F$ can be chosen as close as desired to $f$ in $\mathcal{C}^k(\bar{D}_0\times P)$.
\end{proposition}

\proof
Let $m$ and $M$ denote the complex dimensions of $X$ and $Z$, respectively. We may assume that $\bar{D}_1$ is a domain in $\CC^m$, that $P_0$ is a polydisc, and that $f(\bar{D}_{0,1}\times P_0)$ and $f_0(\bar{D}_1)$ are contained in the same special coordinate chart $\Omega \subset Z$, $\Omega \cong \bar{D}_1\times V$. Then we have
\begin{eqnarray*}
f(x,t) & = & (x,\tilde{f}(x,t)),\ \quad x\in \bar{D}_{0,1},\ t\in P_0, \\
f_0(x) & = & (x,\tilde{f}_0(x,t)), \quad x\in \bar{D}_1,\quad t\in P_0.
\end{eqnarray*}
By Taylor expansion in the variable $t\in \CC ^N$ we have
$$\tilde{f}(x,t)=\tilde{f}_0(x)+\sum _{j=1}^Ng_j(x)t_j+\sum _{j,k=1}^Nh_{j,k}(x,t)t_jt_k, \quad x\in \bar{D}_{0,1}, t\in P_0$$
for some maps $h_{j,k}:\bar{D}_{0,1}\times P_0\to \CC ^{M-m}$ of class $\AA ^k$. Set
$$L_x(t)=\sum _{j=1}^Ng_j(x)\, t_j.$$
Then $L_x=\d _t|_{t=0}\tilde{f}(x,t):\CC ^N\to \CC ^{M-m}$ is a linear map for all $x\in \bar{D}_{0,1}$.
Chose a polydisc $0\in P'\Subset P_0$ and apply Mergelyan's theorem to approximate maps $h_{j,k}$ in $\mathcal{C}^k(\bar{D}_{0,1}\times \bar{P}')$ by entire maps $\hat{h}_{j,k}:\CC ^m\times \CC ^N\to \CC ^{M-m}$. Similarly we approximate maps $g_j$ in $\mathcal{C}^k(\bar{D}_{0,1})$ by entire maps $\hat{g}_j:\CC ^m\to \CC ^{M-m}$. Thus for every $x\in \CC ^m$ we get a linear map
$$\hat{L}_x:\CC ^N\to \CC ^{M-m}, \quad \hat{L}_x(t)=\sum _{j=1}^N\hat{g}_j(x)t_j.$$
We may assume that this map is surjective for every $x\in \bar{D}_1$. Indeed, the set $Q$ of all non-surjective linear maps in $\textrm{Hom}(\CC ^N,\CC ^{M-m})$ has codimension $N-(M-m)+1$ (this follows for example from \cite[p.\ 324, Lemma 7.9.2]{For3}). Applying the basic transversality theorem (see \cite[p.\ 317, Theorem 7.8.5]{For3}) to the map $x\to \hat{L}_x$ gives the desired result provided that $N-(M-m)+1>m$, which is true by assumption. Transversality in this case means that the image of the map does not intersect $Q$.

Consider the map 
$$\tilde{G}(x,t)=\tilde{f}_0(x)+\hat{L}_x(t)+\sum _{j,k=1}^N\hat{h}_{j,k}(x,t)t_jt_k, \quad x\in \bar{D}_1, t\in \CC ^N.$$
Choose an open set $0\in P_1\subset P'$ such that the image $\tilde{G}(\bar{D}_1\times P_1)$ is contained in $V$. Then the map $G:\bar{D}_1\times P_1\to Z$ defined by $G(x,t)=(x,\tilde{G}(x,t))$ is a spray of $\pi $-sections of class $\AA ^k$ with core $f_0|_{\bar{D}_1}$ which is dominating on $\bar{D}_1$ and close to $f$ in $\mathcal{C}^k(\bar{D}_{0,1}\times P_1)$ (we may have to shrink $P_1$). We can now apply Proposition \ref{gs} to glue these two sprays together into a final spray $F:\bar{D}\times P\to Z$ which has all the desired properties.
\endproof

Notice that the condition $N\ge \dim _{\CC }Z$ is only needed for the dominability of the final spray.

\proof[Proof of Theorem \ref{eost}]
Choose a plurisubharmonic function $\phi :V\to [0,\infty )$ of class $\mathcal{C}^l$, defined on an open neighborhood $V$ of $\bar{D}$, which is strictly plurisubharmonic outside $E$, it satisfies $E=\set{x\in V}{\phi (x)=0}$, $D=\set{x\in V}{\phi (x)<c}$ for some $c>0$, and has only isolated critical values different from $c$.

Condition $\mathcal{E}$ (see Sec.\ \ref{sec:sprays}) provides us with an open neighborhood $U\subset D$ of $E$ and a holomorphic spray of $\pi $-sections $f:U\times P\to Z$ ($0\in P\subset \CC ^N$) with core $f_0|_U$ which is dominating on $U\setminus E$ and fixes $E$. We may assume that $N\ge \dim _{\CC }Z$, otherwise we just add a few variables and make the spray independent of them. Choose a number $c_0>0$ which is a regular value of $\phi$ such that $D_0=\set{x\in V}{\phi (x)\le c_0}\subset U$.

By \cite[Corollary 2.8]{Hen&Lei} there exists a finite sequence of 1-convex domains $D_1\subset D_2\subset \ldots \subset D_m=D$ such that for each $j=1,\ldots ,m$ we have $D_j=D_{j-1}\cup B_j$ where $(D_{j-1},B_j)$ is a 1-convex Cartan pair of class $\mathcal{C}^l$. Moreover, we can ensure that for every $j=1,\ldots ,m$ the set $\bar{B}_j$ in contained in some coordinate chart on $X$ and the set $f_0(\bar{B}_j)$ is contained in some special coordinate chart on $Z$.

Using Proposition \ref{eacp} applied inductively to $(D_{j-1},B_j)$, $j=1,\ldots ,m$, we can extend the spray of $\pi $-sections $f$ from $\bar{D}_0$ to $\bar{D}_1$, from $\bar{D}_1$ to $\bar{D}_2$, and so on until we reach $\bar{D}_m=\bar{D}$.
\endproof

As an immediate corollary to Theorem \ref{eost} we get the version without boundary.

\begin{corollary}\label{eosc}
Let $X$ be a 1-convex manifold with the exceptional set $E$, $D\subset X$ a relatively compact domain containing $E$, and $\pi :Z\to X$ a holomorphic submersion which satisfies Condition $\mathcal{E}$. Given a holomorphic $\pi $-section $f_0:X\to Z$ there exists a holomorphic spray of $\pi $-sections $f:D\times P\to Z$ with core $f_0|_{D}$ which is dominating on $D\setminus E$ and fixes $E$.
\end{corollary}

\section{Approximation of sections}
\label{sec:approximation}

We are now ready to prove the main results stated in Sec.\ \ref{sec:main} and to obtain some generalizations. We follow the methods developed in \cite{Dri&For1} for the case of Stein domains, combined with the technique of sprays over $1$-convex domains developed in the previous sections. We begin with a local Mergelyan-type approximation result analogous to \cite[Theorem 5.1]{Dri&For1}. 

\begin{theorem}\label{la}
Let $D$ be a relatively compact 1-convex domain with strongly pseudoconvex boundary of class $\mathcal{C}^l$ $(l\ge 2)$ and exceptional set $E$ in a complex manifold $X$, and let $\pi :Z\to X$ be a holomorphic submersion which satisfies Condition $\mathcal{E}$ (see Sec.\ \ref{sec:sprays}). Every $\pi $-section $f:\bar{D}\to Z$ of class $\AA ^k$ ($k\le l$) can be approximated in $\mathcal{C}^k(\bar{D})$ by $\pi $-sections which are holomorphic in open neighborhoods of $\bar{D}$ and agree with $f$ on $E$. Moreover, the approximating $\pi $-sections can be chosen homotopic to $f$ on $\bar{D}$ through a homotopy of $\pi $-sections of class $\AA ^k(D)$ which agree with $f$ on $E$.
\end{theorem}

\proof 
The proof is practically the same as in \cite{Dri&For1} and we include it for the sake of completeness.
Let $m$ and $M$ denote the complex dimensions of $X$ and $Z$ respectively. Denote by $pr _1:\CC ^m\times \CC ^{M-m}\to \CC ^m$, $pr _2:\CC ^m\times \CC ^{M-m}\to \CC ^{M-m}$ the coordinate projections and by $B\subset \CC ^m$, $B'\subset \CC ^{M-m}$ the unit balls. Since $\pi $ is a holomorphic submersion, there exist for each point $z\in Z$ open neighborhoods $z\in W\subset Z$, $\pi (z)\in V=\pi(W)\in X$ and biholomorphic maps $\Phi :W\to B\times B'$, $\phi :V\to B$ such that
$$\Phi(z)=\Big(\phi (\pi (z)),\phi '(z)\Big)\in B\times B', \quad z\in W,$$
where $\phi '=pr _2\circ \Phi $. Such $(W,\Phi )$ is called a special coordinate chart on $Z$.

Fix a $\pi$-section $f:\bar{D}\to Z$ of class $\AA ^k$. By Narasimhan's lemma on local convexification we can find finitely many special coordinate charts $(W_j,\Phi _j)$ on $Z$, $j=1,\ldots ,J$, with $V_j=\pi (W_j)$ and $\Phi _j=(\phi _1\circ \pi ,\phi _j')$ as above, such that $\set{V_j}{j=1,\ldots ,J}$ is a covering for $\d D$ which does not intersect $E$ and we have
\begin{enumerate}
\item $\phi _j(\d D\cap V_j)$ is a strongly convex hypersurface in $B$,
\item $f(\bar{D}\cap V_j)\subset W_j$ and $\overline{\phi _j'(f(D\cap V_j))}\subset B'$.
\end{enumerate}
Choose a number $0<r<1$ such that the sets $U_j:=\phi _j^{-1}(rB)$, $j=1,\ldots ,J$, still cover $\d D$.

By induction we shall construct a sequence $D=D_0\subset D_1\subset \ldots \subset D_J\Subset X$ of 1-convex domains with strongly pseudoconvex boundary of class $\mathcal{C}^l$ (and the same exceptional set) and $\pi $-sections $f_j:\bar{D}_j\to Z$ of class $\AA ^k$, with $f_0=f$, such that $f_j$ is close to $f_{j-1}$ in $\mathcal{C}^k(D_{j-1})$ and agrees with $f$ on $E$ for every $j=1,\ldots ,J$. To keep the induction going we will insure that the property (ii) will remain valid for $(D_j,f_j)$ in place of $(D,f)$. Due to this fact the domain $D_j$ will in general depend on $\mathcal{C}^0(\bar{D}_{j-1})$ distance between $f_j$ and $f_{j-1}$ and it will be chosen so that
$$D_{j-1}\subset D_j\subset D_{j-1}\cup V_j, \quad \d D_{j-1}\cap U_j\subset D_j.$$
Thus the final domain $D_J$ will contain $\bar{D}$ in its interior and the final section $f_J$ will approximate $f$ as close as desired in $\mathcal{C}^k(\bar{D})$ and agree with $f$ on $E$.

\begin{figure}[h]

\psset{unit=0.7cm}
\begin{pspicture}(-4.1,-4.1)(4.1,4.1)
\SpecialCoor

\pspolygon[fillstyle=solid,fillcolor=lightgray,linecolor=lightgray](3;117.8)(3.5;135)(3.5;-135)(3;-117.8)

\psecurve[fillstyle=solid,fillcolor=lightgray](4;129)(3;117.8)(1.5,0)(3;-117.8)(4;-129)
\pscurve(4.5;130)(1,0)(4.5;-130)
\psecurve(4;135)(3;117.8)(0.5,0)(3;-117.8)(4;-135)

\psecurve[fillstyle=solid,fillcolor=lightgray](2.5;120)(3;117.8)(3.5;135)(3;152.2)
\psarc[fillstyle=solid,fillcolor=lightgray](0,0){3.5}{135}{-135}
\psecurve[fillstyle=solid,fillcolor=lightgray](2.5;-120)(3;-117.8)(3.5;-135)(3;-152.2)

\pscircle[linestyle=dashed](0,0){4}
\pscircle[linestyle=dashed](0,0){2}

\rput(1;180){$\d D$}
\psline[arrows=->](0.6;180)(0.4,0)
\rput(3;180){$D_1'$}
\rput(3;0){$\d D_1$}
\psline[arrows=->](2.5;0)(1.1,0)
\rput(3.5;90){$V_1$}
\rput(3;45){$U_1$}
\psline[arrows=->](2.7;45)(2.1;45)
\end{pspicture}

\caption{The domains $D_1'$ and $D_1$ (modified from \cite[Fig.\ 1]{Dri&For1}).}
\label{figD}
\end{figure}
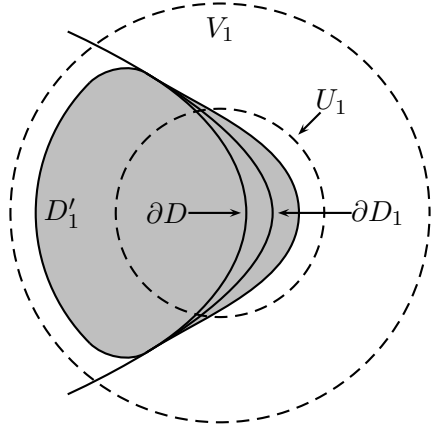

All the steps in the induction are the same, so we only need to explain how to get $(D_1,f_1)$ from $(D,f)$. We first find a domain $D_1'\subset V_1$ with $\mathcal{C}^l$ boundary which is a convex bump on $D$ such that $\bar{D}\cap \bar{U}_1\subset D_1'$. Choose a smooth function $\chi :\CC ^m\to [0,1]$ with compact support contained in $B$ such that $\chi =1$ on $rB$. Let $\tau :B\to \RR $ be a strictly convex defining function of class $\mathcal{C}^l$ for the hypersurface $\phi _j(\d D\cap V_j)$. We may assume it has no critical points. Choose $r_0\in (r,1)$ close to $1$ such that the hypersurface $\{\tau =0\}=\phi _j(\d D\cap V_j)$ intersects the sphere $\d (r_0B)$ transversally. Choose $\delta >0$ small enough so that $\tau -\delta \chi $ is still a strictly convex function without critical points. By rounding off the corners of the domain $r_0B\cap \{\tau -\delta \chi <0\}$ (i.e. intersection of $\{\tau -\delta \chi =0\}$ with $\d (r_0B)$) we get a domain $D_1''$ such that the domain $D_1'=\phi ^{-1}(D_1'')$ has the desired properties (see Figure \ref{figD}).

By Theorem \ref{eost} there exists a spray of $\pi $-sections $F:\bar{D}\times P\to Z$ ($0\in P\subset \CC ^N$) of class $\AA ^k$ with core $f|_{\bar{D}}$ which is dominating on $\bar{D}\setminus E$ and fixes $E$. Shrinking $P$ around $0$ if necessary we may assume that property (ii) remains valid for $F_t$ in place of $f$ for all $t\in P$. By using the special coordinate chart $(W_1,\Phi _1)$ we can find an open set $\Omega \subset V_1$ containing $\bar{D}_1'\cap \bar{D}$ and a holomorphic spray of $\pi $-sections $G:\Omega \times P_0$ ($0\in P_0\subset P$) with image contained in $W_1$ which is close to $F$ in $\mathcal{C}^k(\bar{D}_1'\cap\bar{D}\times P_0)$. This is done in the same way as in the proof of the Proposition \ref{eacp} except that we also need to approximate the first coefficient of the Taylor expansion (since $f=F_0$ is only defined on $\bar{D}$). Since the image of this new coefficient needs to stay inside $W_1$ it may happen that $\Omega $ is only a small neighborhood of $\bar{D}_1'\cap \bar{D}$ and does not contain the whole $\bar{D}_1'$.

We would like to use Proposition \ref{gs} to glue the sprays of $\pi $-sections $F$ and $G$ into a single spray of $\pi $-sections. We can not do this directly since their domains do not form a 1-convex Cartan pair. However, we can still use Lemma \ref{tm} to find a transition map $\gamma :\bar{D}_1'\cap\bar{D}\times P_1\to \CC ^N$ ($0\in P_1\subset P_0$) of class $\AA ^k$ between $F$ and $G$ which is close to $\gamma _0(x,t)=t$ in $\mathcal{C}^k(\bar{D}_1'\cap\bar{D})$ and satisfies $F(x,t)=G(x,\gamma (x,t))$ for $x\in \bar{D}_1'\cap\bar{D}$, $t\in P_1$. Theorem \ref{st} applied to $\gamma $ on the 1-convex Cartan pair $(\bar{D},\bar{D}_1')$ furnishes a set $0\in P_2\subset P_1$ and maps $\alpha :\bar{D}\times P_2\to \CC ^N$, $\beta :\bar{D}_1'\times P_2\to \CC ^N$ of class $\AA ^k$ which are close to $\gamma _0$ in $\mathcal{C}^k(\bar{D})$, $\mathcal{C}^k(\bar{D}_1')$, respectively, and satisfy $\gamma (x,\alpha (x,t))=\beta (x,t)$ for $x\in \bar{D}_1'\cap\bar{D}$, $t\in P_2$. The sprays of $\pi $-sections $F(x,\alpha (x,t))$ and $G(x,\beta (x,t))$ then amalgamate into a single spray of $\pi $-sections defined on $\bar{D}\cup (\bar{D}_1'\cap \Omega )$ which is close to $F$ in $\mathcal{C}^k(\bar{D}\times P_2)$ and fixes $E$. The core of this new spray of $\pi $-sections, which we denote by $f_1$, is a $\pi $-section of class $\AA ^k$ which is close to $f$ in $\mathcal{C}^k(\bar{D})$ and agrees with $f$ on $E$.

It remains to restrict $f_1$ to a suitably chosen 1-convex domain $D_1\Subset X$ contained in $\bar{D}\cup (\bar{D}_1'\cap \Omega )$ which satisfies all the required properties. We choose $D_1$ such that it agrees with $D$ outside $V_1$ and
$$D_1\cap V_1=\phi _1^{-1}(\{\tau -\epsilon \chi <0\})$$
for a $0<\epsilon <\delta $. By choosing $\epsilon >0$ small enough we ensure that properties (i) and (ii) are satisfied for $(D_1,f_1)$. This concludes the induction step.

Notice that despite the fact that we were not able to use Proposition \ref{gs} directly we were still able to use the same technique to obtain the approximating $\pi $-section. Hence we also get that the approximating $\pi $-section is homotopic to $f$ on $\bar{D}$ through a homotopy of $\pi $-sections of class $\AA ^k(D)$ which agree with $f$ on $E$.
\endproof

The domains of the approximating sections in Theorem \ref{la} must in general shrink to $\bar{D}$. If we wish to get approximation result by global sections (as in Theorem \ref{opr}) we need to make additional assumptions. To avoid topological obstructions we assume that the initial section can be extended continuously to a global one. More importantly we will need to assume that the submersion $\pi $ is in fact a fiber bundle with Oka fiber (see Definition \ref{om}). This will compensate for the fact that the initial $\pi $-section (core of the spray of $\pi $-sections) is not holomorphic on the domain where we want to extend the spray of $\pi $-sections.

The next result is a direct corollary of Theorem \ref{la} and the main result in \cite{Pre}. However, it also follows from our main Theorem \ref{opr} (independently of \cite{Pre}) by a simple induction argument, which we sketch in the proof.

\begin{corollary}\label{opc}
Let $X$ be a 1-convex manifold with exceptional set $E$ and $D\Subset X$ be a 1-convex domain with $\mathcal{C}^l$ ($l\ge 2$) boundary and $\OO (X)$-convex closure such that $E\subset D$. Let $\pi :Z\to X$ be a holomorphic submersion which satisfies Condition $\mathcal{E}$ such that $\pi :Z|_{X\setminus \bar{D}}\to X\setminus \bar{D}$ is a holomorphic fiber bundle with Oka fiber. For every continuous $\pi $-section $f_0:X\to Z$ of class $\AA ^k$ ($k\le l$) on $\bar{D}$ there exists a homotopy $f_t:X\to Z$ of continuous $\pi $-sections of class $\AA ^k$ on $\bar{D}$ such that $f_t$ is close to $f_0$ in $\mathcal{C}^k(\bar{D})$ and agrees with $f_0$ on $E$ for every $t\in [0,1]$, and $f_1$ is holomorphic on $X$.
\end{corollary}

\proof[Proof (sketch)] For simplicity we assume that $\pi :Z\to X$ is a holomorphic fiber bundle with Oka fiber over the whole manifold $X$. Choose a smooth exhaustion function $\phi :X\to \RR $ which is strictly plurisubharmonic Morse function outside $E$ such that $\phi >0$ on $\bar{D}$. Choose a sequence of numbers $0<c_1<c_2<\ldots $ which are not critical values of $\phi $ such that $\lim _{j\to \infty }c_j=\infty $. Set $D_0:=D$ and $D_j:=\set{x\in X}{\phi (x)<c_j}$ for $j\in \NN $. We divide the parameter interval $[0,1]$ of the homotopy into subintervals $I_j=[t_j,t_{j+1}]$ with $t_j=1-2^{-j}$, $j=0,1,2,\ldots $. Fix an $\epsilon >0$. We shall construct a homotopy of continuous $\pi $-sections $f_t:X\to Z$ ($t\in [0,1)$) such that for every $j=0,1,2,\ldots $ and $t\in I_j$ the $\pi $-section $f_t$ is of class $\AA ^k$ on $\bar{D}_j$, agrees with $f_0$ on $E$, and satisfies
$$\sup _{t\in I_j}\| f_t-f_{t_j}\| _{\mathcal{C}^k(\bar{D}_j)}<2^{-j-1}\epsilon .$$
The limit section $f_1:=\lim _{t\to 1}f_t:X\to Z$ is then holomorphic on $X$ and $\epsilon $ close to $f_0$ in $\mathcal{C}^k(\bar{D})$.

We assume inductively that the homotopy has already been constructed for $t\in [0,t_j]$ and explain how to construct it for $t\in I_j$. The continuous $\pi $-section $f_{t_j}:X\to Z$ is of class $\AA ^k$ on $\bar{D}_j$. By using Theorem \ref{opr} for a pair of domains $D_j\Subset D_{j+2}\Subset X$ we find a homotopy of continuous $\pi $-sections $\tilde{f}_t:\bar{D}_{j+2}\to Z$ for $t\in I_j$ which has all the desired properties, except that it is not yet defined on the whole $X$. It is easy to correct this homotopy to a global one, keeping it unchanged on $\bar{D}_{j+1}$ and at $t=t_j$. Just take a smooth function $\chi :X\to [0,1]$ with compact support contained in $D_{j+2}$ such that $\chi =1$ on $\bar{D}_{j+1}$ and define
$$f_t(x):=\tilde{f}_{t\chi (x)+t_j(1-\chi (x))}(x).$$
\endproof

Corollary \ref{opc} is an 'open' version of Theorem \ref{opr} where the closed base domain $\bar{D}'$, over which the fiber bundle is defined, is replaced by a manifold $X$ (without boundary).

We now proceed to the proof of Theorem \ref{opr}. The ideas are similar as in the proof of the existence of sprays of sections. We `thicken' the initial section into a spray of sections and then extend it to bigger and bigger domains in the base manifold. There will be two major steps in the proof: we will need to explain how to extend across a convex bump and across a critical point. The main difference is that before the initial section (the core of the spray of sections) was already defined and holomorphic on the bigger domain where we wanted to extend the spray of sections. Now this will not be the case, hence we will not be able to keep it fixed.

\begin{proposition}\label{eacb}
Let $(D_0,D_1)$ be a 1-convex Cartan pair of class $\mathcal{C}^l$ ($l\ge 2$) in a complex manifold $X$, where $D_1$ is a convex bump on $D_0$. Set $D_{0,1}=D_0\cap D_1$ and $D=D_0\cup D_1$ and denote by $E$ the exceptional set of $D_0$. Let $\pi:Z\to \bar{D}$ be either a fiber bundle of class $\AA ^k$ ($k\le l$) or the restriction of a holomorphic submersion $\tilde{Z}\to X$, such that $\pi :Z|_{D_0}\to D_0$ satisfies Condition $\mathcal{E}$ and $\pi :Z_{\bar{D}_1}\to \bar{D}_1$ is a trivial fiber bundle with Oka fiber. For every continuous $\pi $-section $f_0:\bar{D}\to Z$ of class $\AA ^k$ on $\bar{D}_0$ there exists a homotopy $f_t:\bar{D}\to Z$ of continuous $\pi $-sections of class $\AA ^k$ on $\bar{D}_0$ such that $f_t$ is close to $f_0$ in $\mathcal{C}^k(\bar{D}_0)$ and agrees with $f_0$ on $E$ for every $t\in [0,1]$, and $f_1$ is of class $\AA ^k$ on $\bar{D}$.
\end{proposition}

\proof
By Theorem \ref{eost} there exists a spray of $\pi $-sections $F:\bar{D}_0\times P\to Z$ of class $\AA ^k$ with core $f$ which is dominating on $\bar{D}_0\setminus E$ and fixes $E$. Using the trivialization $Z|_{\bar{D}_1}\cong \bar{D}_1\times Y$ ($Y$ is the fiber) we can write
$$F(x,t)=(x,\tilde{F}(x,t))\in \bar{D}_{0,1}\times Y, \quad x\in \bar{D}_{0,1}, t\in P.$$
By the definition of a convex bump we may assume that $D_{0,1}$ and $D_1$ are strongly convex domains in $\CC ^m$. Choose a polydisc $0\in P_0\Subset P$. For $0<\lambda <1$ the map $\tilde{F}_{\lambda }(x,t)=\tilde{F}(\lambda x,t)$ is holomorphic in $\frac{1}{\lambda }D_{0,1}\times P\supset \bar{D}_{0,1}\times \bar{P}_0$. If $\lambda $ is chosen close enough to $1$ than $\tilde{F}_{\lambda }$ is as close as desired to $\tilde{F}$ in $\mathcal{C}^k(\bar{D}_{0,1}\times P)$. Since $Y$ is an Oka manifold we can approximate $\tilde{F}_{\lambda }$ in $\mathcal{C}^k(\bar{D}_{0,1}\times \bar{P}_0)$ by an entire map $\tilde{G}$ and set
$$G(x,t)=(x,\tilde{G}(x,t))\in \bar{D}_1\times Y, \quad x\in \bar{D}_1, t\in P.$$
If the above approximations are close enough we can use Proposition \ref{gs} to glue the sprays of $\pi $-sections $F$ and $G$ into a single spray of $\pi $-sections $H:\bar{D}\times P_1\to Z$ ($0\in P_1\subset P_0$) of class $\AA ^k$ which is close to $F$ in $\mathcal{C}^k(\bar{D}_0\times P_1)$ and fixes $E$. The core $H_0=:f_1$ of this spray of $\pi $-sections is than of class $\AA ^k$ on $\bar{D}$, close to $f_0$ in $\mathcal{C}^k(\bar{D}_0)$, and agrees with $f_0$ on $E$. Moreover, $f_1$ is homotopic to $f_0$ on $\bar{D}_0$ through a homotopy with desired properties. It is not hard to see that this homotopy extends to $\bar{D}$. For details see \cite{For&Pre}.
\endproof

\begin{proposition}\label{eah}
Let $Z$ be a complex manifold, $X$ a 1-convex manifold with the exceptional set $E$, and $\pi :Z\to X$ a holomorphic submersion which satisfies Condition $\mathcal{E}$. Let $D\Subset X$ be a 1-convex domain with $E\subset D$ and let $M\subset X$ be a $\mathcal{C}^1$ totally real submanifold such that $\bar{D}$ and $\bar{D}\cup M$ are compact $\OO (X)$-convex subsets. Given an open neighborhood $U\subset X$ of $\bar{D}$, a continuous $\pi $-section $f:U\cup M\to Z$ which is holomorphic on $U$, an $\epsilon >0$, and a $k\ge 0$, there exists an open neighborhood $V\subset X$ of $\bar{D}\cup M$ and holomorphic $\pi $-section $g:V\to Z$ which agrees with $f$ on $E$ and is $\epsilon $ close to $f$ in $\mathcal{C}^k(\bar{D})$ and $\mathcal{C}(M)$. Moreover, $g$ can be chosen homotopic to $f$ through a homotopy of $\pi $-sections which are holomorphic in neighborhoods of $\bar{D}$ and continuous in neighborhoods of $M$.
\end{proposition}

\proof
Choose an open neighborhood $U_0\Subset U$ of $\bar{D}$. By Corollary \ref{eosc} there exists a holomorphic spray of $\pi $-sections $F:U_0\times P\to Z$ with core $f|_{U_0}$ which is dominating on $U_0\setminus E$ and fixes $E$. Shrinking $U_0$ and $P$ if necessary we can extend the spray $F$ continuously to $M\times P$.

By a standard result (see for example \cite[p. 72, Corollary 3.5.2]{For3}) the set $B:=M\setminus D\subset X$ has a basis of Stein neighborhoods. On the other hand the set $\bar{D}\cup M$ has a basis of 1-convex neighborhoods. For $D$ Stein (i.e., $E=\emptyset $) this was proved in \cite[Theorem 3.1]{For} by gluing certain plurisubharmonic functions (see also \cite[p. 78, Theorem 3.7.1]{For3}). The same proof applies in our 1-convex case, since the gluing takes place only in a neighborhood of the set $M\cap \d D$, i.e., away from the exceptional set $E$.

Fix a Stein neighborhood $\Omega \Subset X$ of $B$ with smooth strictly pseudoconvex boundary such that $\bar{\Omega }\cap E=\emptyset $. Moreover, fix a 1-convex neighborhood $D_0\Subset U_0$ of $\bar{D}$ with smooth strictly pseudoconvex boundary. Set $K:=\bar{D}_0\cap \bar{\Omega }$ and $S:=K\cup B$. Now we can use \cite[Theorem 3.2]{For} to approximate the spray of $\pi $-sections $F$, which is holomorphic over a neighborhood $U_0$ of $K$ and continuous over $B$, by a spray of $\pi $-sections $G:\Omega _0\times P_0\to Z$ ($0\in P_0\Subset P$), which is holomorphic over some neighborhood $\Omega _0$ of $S$ and close to $F$ in $\mathcal{C}^k(K\times P_0)$ and $\mathcal{C}(B\times P_0)$. If the approximation is close enough then the core $G_0$ is homotopic to $f=F_0$ through a homotopy of $\pi $-sections which are holomorphic in a neighborhood of $K$ and continuous in a neighborhood of $B$.

It remains to glue the sprays of $\pi $-sections $F$ and $G$ into a single spray of $\pi $-sections. To do this choose a 1-convex neighborhood $V$ of $\bar{D}\cup M$ with smooth strictly pseudoconvex boundary such that $V\Subset D_0\cup (\Omega \cap \Omega _0)$. By smoothing the corners of $V\cap D_0$ we get a 1-convex domain $A$ and by smoothing the corners of $V\cap \Omega $ we get a Stein domain $B$. This two domains form a 1-convex Cartan pair $(A,B)$ of class $\mathcal{C}^{\infty }$ with $A\cup B=V$, $D\Subset A\Subset U_0$, $B\Subset \Omega _0$, and $A\cap B\subset K$. Using Proposition \ref{gs} we can now glue the sprays of $\pi $-sections $F|_{\bar{A}\times P_0}$ and $G|_{\bar{B}\times P_0}$ into a single spray of $\pi $-sections $H:\bar{V}\times P_1\to Z$ ($0\in P_1\subset P_0$). Its core $g:=G_0|_V$ is then a holomorphic $\pi $-section which satisfies all the conclusions in the proposition.
\endproof

\proof[Proof of Theorem \ref{opr}]
First we show that we may assume $f_0$ to be holomorphic in some open neighborhood $U\subset D'$ of $\bar{D}$. If this is not the case, we can use Theorem \ref{la} to find a $\pi $-section $g$ holomorphic in an open neighborhood $U_0\subset D'$ of $\bar{D}$ which is close to $f$ in $\mathcal{C}^k(\bar{D})$ and agrees with $f_0$ on $E$. Moreover, if the approximation is close enough and $U_0$ is chosen small enough then $g$ is homotopic to $f_0$ on $U_0$ through a homotopy $g_t:U_0\to Z$ of continuous $\pi $-sections which are of class $\AA ^k$ on $\bar{D}$ and agree with $f_0$ on $E$. Choose a smooth function $\chi :X\to [0,1]$ with compact support contained in $U_0$ such that $\chi =1$ in an open neighborhood $U$ of $\bar{D}$. Define $\tilde{g}_t(x):=g_{t\chi (x)}(x)$. This is a homotopy of global continuous $\pi $-sections which are of class $\AA ^k$ on $\bar{D}$ and agree with $f_0$ on $E$, such that $\tilde{g}_1$ is holomorphic on $U\subset \bar{D}$. Thus we can replace $f_0$ by $\tilde{g}_1$.

The assumptions guarantee that there exists a strictly plurisubharmonic function $\phi $ of class $\mathcal{C}^l$ defined on an open neighborhood $V\subset X$ of $\bar{D}'\setminus D$ which satisfies $D\cap V\subset \set{x\in V}{\phi (x)<0}\Subset U$ and $D'\cap V=\set{x\in V}{\phi (x)<1}$. Moreover, we may assume that $\phi $ is a Morse function with {\em nice critical points} $p_1,\ldots ,p_J\in V$ (in the sense of Def.\ 3.9.2 in \cite[p.\ 89]{For3}) such that $0<\phi (p_1)<\ldots <\phi (p_J)<1$. Choose a sequence of numbers $0<c_1<\ldots <c_{2J}<1$ such that $c_{2j-1}<\phi (p_j)<c_{2j}$ for every $j=1,\ldots ,J$ and define $c_0:=0$, $c_{2J+1}=1$. Denote $D_j=D\cup \set{x\in V}{\phi (x)<c_j}$, $j=0,\ldots ,2J+1$. Divide the parameter interval $[0,1]$ of the homotopy into subintervals $I_j=[t_j,t_{j+1}]$, $j=0,\ldots ,2J$, with $t_j=\frac{j}{2J+1}$, $j=0,\ldots ,2J+1$.

We shall construct a homotopy $f_t:\bar{D}'\to Z$ of $\pi $-sections such that for every $j=0,\ldots ,2J$ and $t\in I_j$ the $\pi $-section $f_t$ is of class $\AA ^k$ on $\bar{D}_j$, agrees with $f_0$ on $E$, and is $\frac{\epsilon }{2J+1}$ close to $f_{t_j}$ in $\mathcal{C}^k(\bar{D}_j)$. Thus $f_t$ will be $\epsilon $ close to $f_0$ in $\mathcal{C}^k(\bar{D})$ for every $t\in [0,1]$.

We assume inductively that the homotopy has already been constructed for $t\in [0,t_j]$ and explain how to construct it for $t\in I_j=[t_j,t_{j+1}]$. We consider two cases.

\medskip
\noindent {\it The non-critical case:} If $j$ is even then $\phi $ has no critical points in $\bar{D}_{j+1}\setminus D_j$. By the bumping lemma (see for example \cite[Lemma 2.2]{Hen&Lei} or \cite[Lemma 5.10.3, p.\ 218]{For3}) there exists a finite sequence 
$$D_j=\Omega _0\subset \Omega _1\subset \ldots \subset \Omega _m=D_{j+1}$$
(where $m$ depends on $j$) of 1-convex domains with strongly pseudoconvex boundary of class $\mathcal{C}^l$ such that for each $i=1,\ldots ,m$ we have $\Omega_i=\Omega _{i-1}\cup B_i$, where $B_i$ is a convex bump on $\Omega _{i-1}$. Moreover, we can ensure that $\pi :Z|_{\bar{B}_i}\to \bar{B}_i$ is a trivial fiber bundle for every $i$. Using Proposition \ref{eacb} applied inductively to Cartan pairs $(\Omega_{i-1},B_i)$, $i=1,\ldots ,m$ we get the desired homotopy on the interval $[t_j,t_{j+1}]$.

\medskip
\noindent {\it The critical case:} If $j$ is odd then $\phi $ has a unique critical point $p\in D_{j+1}\setminus \bar{D}_j$. Denote by $M$ the stable manifold through $p$ of the gradient flow of $\phi $. For $\delta >0$ small enough \cite[p.\ 96, Proposition 3.10.4]{For3} (see also \cite[Lemma 3.1]{Dri&For2}) gives a strictly plurisubharmonic function $\psi $ on $\set{x\in V}{\phi (x)<\phi (p)+3\delta }$ of class $\mathcal{C}^l$ satisfying
\begin{enumerate}
\item $\{\phi \le c_j \}\cup M\subset \{\psi \le 0\}\subset \{\phi \le \phi (p)-\delta \}\cup M$,
\item $\{\phi \le \phi (p)+\delta \}\subset \{\psi \le 2\delta \}\subset \{\phi <\phi (p)+3\delta \}$,
\item $\psi $ has no critical values in $(0,3\delta )$,
\item $\psi=\phi +s$ outside some neighborhood of $p$ for some $s>\delta $.
\end{enumerate}

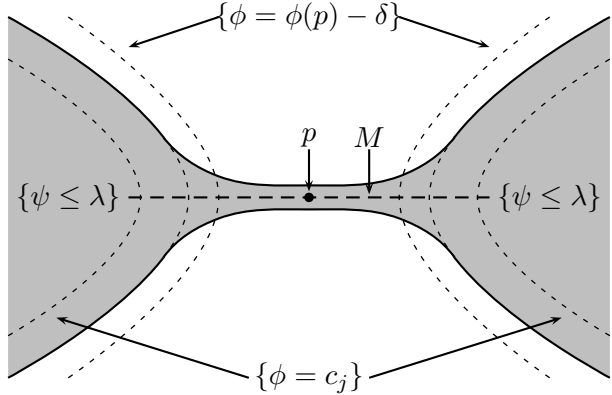
\begin{figure}[h]

\psset{unit=0.8cm}
\begin{pspicture}(-5.1,-3.1)(5.1,3.1)

\pscurve[fillstyle=solid,fillcolor=lightgray](-5,3)(-2,0)(-5,-3)
\pscurve[fillstyle=solid,fillcolor=lightgray](5,3)(2,0)(5,-3)

\pspolygon[fillstyle=solid,fillcolor=lightgray,linestyle=none](-2.5,0.95)(-2.5,-0.95)(2.5,-0.95)(2.5,0.95)

\psline[linestyle=dashed](-3,0)(3,0)

\psecurve[linestyle=dashed,dash=2pt 3pt,linewidth=0.5pt](-5,3)(-2.42,1)(-2,0)(-2.42,-1)(-5,-3)
\psecurve[linestyle=dashed,dash=2pt 3pt,linewidth=0.5pt](5,3)(2.42,1)(2,0)(2.42,-1)(5,-3)

\psecurve[fillstyle=solid,fillcolor=white](-2.8,2)(-2.42,1)(-1.5,0.35)(0,0.2)(1.5,0.35)(2.42,1)(2.8,2)
\psecurve[fillstyle=solid,fillcolor=white](-2.7,-2)(-2.42,-1)(-1.5,-0.35)(0,-0.2)(1.5,-0.35)(2.42,-1)(2.7,-2)

\pscurve[linestyle=dashed,dash=2pt 3pt,linewidth=0.5pt](-4,3)(-1.5,0)(-4,-3)
\pscurve[linestyle=dashed,dash=2pt 3pt,linewidth=0.5pt](4,3)(1.5,0)(4,-3)
\pscurve[linestyle=dashed,dash=2pt 3pt,linewidth=0.5pt](-5,2.3)(-2.8,0)(-5,-2.3)
\pscurve[linestyle=dashed,dash=2pt 3pt,linewidth=0.5pt](5,2.3)(2.8,0)(5,-2.3)

\psdot(0,0)

\rput(0,1){$p$}
\psline[arrows=->](0,0.8)(0,0.1)
\rput(1,1){$M$}
\psline[arrows=->](1,0.8)(1,0.05)
\rput(0,3){$\{\phi =\phi(p)-\delta \}$}
\psline[arrows=->](-1.5,3)(-3,2.5)
\psline[arrows=->](1.5,3)(3,2.5)
\rput(0,-3){$\{\phi =c_j\}$}
\psline[arrows=->](-1,-3)(-4.2,-2)
\psline[arrows=->](1,-3)(4.2,-2)
\rput(4,0){$\{\psi \le \lambda\}$}
\rput(-4,0){$\{\psi \le \lambda\}$}

\end{pspicture}

\caption{The typical level set of the function $\psi $ (modified from \cite[p.\ 94, Fig. 3.5]{For3}).}
\label{figP}
\end{figure}

We subdivide the interval $I_j=[t_j,t_{j+1}]$ into four subintervals $J_1=[t_j,s_1]$, $J_2=[s_1,s_2]$, $J_3=[s_2,s_3]$, $J_4=[s_3,t_{j+1}]$ with $t_j<s_1<s_2<s_3<t_{j+1}$ and proceed in four steps.

\medskip
\noindent {\it Step 1:} By the non-critical case we deform the $\pi $-section $f_{t_j}$ to a $\pi $-section $f_{s_1}$ of class $\AA ^k$ on $D\cup \{\phi \le \phi (p)-\frac{\delta }{2}\}$ through a homotopy $f_t:\bar{D}'\to Z$ ($t\in J_1$) of $\pi $-sections that has the desired properties.

\medskip
\noindent {\it Step 2:} The set $M\setminus (D\cup \{\phi <\phi (p)-\delta \})$ is a $d$-dimensional totally real ball $B\subset X$ attached with its boundary $(d-1)$-sphere $\d B$ to $\{\phi =\phi (p)-\delta \}$. (Here $d$ is the Morse index of the critical point $p$.) We now use Proposition \ref{eah} to approximate the $\pi $-section $f_{s_1}$ in $\mathcal{C}^k(D\cup \{\phi \le \phi (p)-\delta \})$ by a $\pi $-section $f_{s_2}$ holomorphic in some open neighborhood $W_0$ of $D\cup \{\phi \le \phi (p)-\delta \}\cup B$ which agrees with $f_0$ on $E$. Moreover, if the approximation is close enough and $W_0$ is chosen small enough then $f_{s_2}$ is homotopic to $f_{s_1}$ on $W_0$ through a homotopy $f_t:W_0\to Z$ ($t\in J_2$) of holomorphic $\pi $-sections which agree with $f_0$ on $E$. As in the beginning of the proof we can correct this homotopy to a global one, keeping it unchanged in some open neighborhood $W$ of $D\cup \{\phi \le \phi (p)-\delta \}\cup B$ and at $t=s_1$.

\medskip
\noindent {\it Step 3:} By (i) for a small enough $\lambda >0$ the set $D\cup \{\psi \le \lambda \}$ is contained in $W$ (and contains $\bar{D}_j$). Applying the non-critical case with the function $\psi $ (see (iii)) we deform the $\pi $-section $f_{s_2}$ to a $\pi $-section $f_{s_3}$ of class $\AA ^k$ on $D\cup \{\psi \le 2\delta \}$ through a homotopy $f_t:\bar{D}'\to Z$ ($t\in J_3$) of $\pi $-sections that has the desired properties.

\medskip
\noindent {\it Step 4:} By (ii) the $\pi $-section $f_{s_3}$ is holomorphic in a neighborhood of the compact set $\{\phi \le \phi (p)+\delta \}$. Applying the non-critical case we deform it to a $\pi$-section $f_{t_{j+1}}$ of class $\AA ^k$ on $\bar{D}_{j+1}$ through a homotopy $f_t:\bar{D}'\to Z$ ($t\in J_4$) of $\pi $-sections that has the desired properties.

\medskip
\noindent If the approximations in all four steps were close enough then for every $t\in I_j$ the $\pi $-section $f_t$ is $\frac{\epsilon }{2J+1}$ close to $f_{t_j}$ in $\mathcal{C}^k(\bar{D}_j)$. This concludes the proof of the induction step and thus of the theorem. 

Notice that we only needed $\pi $ to be a fiber bundle with Oka fiber outside of the set $\bar{D}$.
\endproof

By adding  a suitable topological condition on the fiber of $\pi$, which will allow us to extend sections continuously to global ones, we get the following corollary to Theorem \ref{opr}. (Compare with \cite[p.\ 233, Theorem 5.14.1]{For3}.) 

\begin{corollary}\label{ga}
Let $X$ be a complex manifold and let $D\Subset D'\Subset X$ be 1-convex domains with strongly pseudoconvex boundary of class $\mathcal{C}^l$ ($l\ge 2$) and the same exceptional set $E$ such that $\bar{D}$ is $\OO (D')$-convex. Let $\pi :Z\to \bar{D}'$ be a fiber bundle of class $\AA ^k$ ($k\le l$) whose fiber $Y$ is an Oka manifold which satisfies $\pi _q(Y)=0$ for all $q<\dim _{\CC }X$. Assume also that the submersion $\pi :Z|_D\to D$ satisfies Condition $\mathcal{E}$. Then every $\pi $-section $f:\bar{D}\to Z$ of class $\AA ^k$ can be approximated in $\mathcal{C}^k(\bar{D})$ by $\pi $-sections $g:\bar{D'}\to Z$ of class $\AA ^k$ which agree with $f$ on $E$.
\end{corollary}

\proof
The only place in the proof of Theorem \ref{opr} where we needed the initial section to be continuous on $X$ is when crossing a critical point $p$ of index $d$ of the strictly plurisubharmonic Morse function $\phi $ (\textit{step 2} of the \textit{critical case}). At such a point we need to be able to extend a given $\pi $-section $f$, defined on the sublevel set $D\cup \{\phi \le \phi (p)-\frac{\delta }{2}\}$, continuously across a $d$-dimensional totally real ball $B\subset X$ attached with its boundary $(d-1)$-sphere $\d B$ to $\{\phi =\phi (p)-\frac{\delta }{2}\}$, such that $D\cup \{\phi \le \phi (p)-\frac{\delta }{2}\}\cup B$ is a strong deformation retract of a sublevel set $D\cup \{\phi \le \phi (p)+\epsilon \}$ for some $\epsilon >0$.

If $\delta $ was chosen small enough then $\bar{B}$ is contained in a small neighborhood of $p$ over which $\pi $ is a trivial fiber bundle with fiber $Y$. Using this trivialization we can identify local $\pi $-sections with maps into $Y$. We see that the required extension exists if and only if the map $f:\d B\to Y$ is null-homotopic. This is certainly true if $\pi _{d-1}(Y)=0$, which holds by the assumption since Morse indices  of a strictly plurisubharmonic function are never bigger than $\dim _{\CC }X$.
\endproof

Similarly we get an analogous corollary to Corollary \ref{opc}.

\medskip
{\it Acknowledgements.} I would like to thank Jasna Prezelj for helpful discussions and Franc Forstneri\v c for introducing me to the problem and for useful suggestions which helped me to improve the paper.

\bibliographystyle{amsplain}

\end{document}